\documentclass[11pt,twoside,a4paper]{article}
\usepackage{makeidx}
\usepackage[dvips]{graphicx}
\usepackage{amscd}
\usepackage{amsmath}
\usepackage{amssymb}
\usepackage{latexsym}
\usepackage[latin1]{inputenc}
\usepackage[T1]{fontenc}   
\usepackage[english]{babel}
\newcommand{\tun}{\begin{picture}(5,0)(-2,-1)
\put(0,0){\circle*{2}}
\end{picture}}

\newcommand{\tdeux}{\begin{picture}(7,7)(0,-1)
\put(3,0){\circle*{2}}
\put(3,0){\line(0,1){5}}
\put(3,5){\circle*{2}}
\end{picture}}

\newcommand{\ttroisun}{\begin{picture}(15,8)(-5,-1)
\put(3,0){\circle*{2}}
\put(-0.65,0){$\vee$}
\put(6,7){\circle*{2}}
\put(0,7){\circle*{2}}
\end{picture}}
\newcommand{\ttroisdeux}{\begin{picture}(5,12)(-2,-1)
\put(0,0){\circle*{2}}
\put(0,0){\line(0,1){5}}
\put(0,5){\circle*{2}}
\put(0,5){\line(0,1){5}}
\put(0,10){\circle*{2}}
\end{picture}}

\newcommand{\tquatreun}{\begin{picture}(15,12)(-5,-1)
\put(3,0){\circle*{2}}
\put(-0.65,0){$\vee$}
\put(6,7){\circle*{2}}
\put(0,7){\circle*{2}}
\put(3,7){\circle*{2}}
\put(3,0){\line(0,1){7}}
\end{picture}}
\newcommand{\tquatredeux}{\begin{picture}(15,18)(-5,-1)
\put(3,0){\circle*{2}}
\put(-0.65,0){$\vee$}
\put(6,7){\circle*{2}}
\put(0,7){\circle*{2}}
\put(0,14){\circle*{2}}
\put(0,7){\line(0,1){7}}
\end{picture}}
\newcommand{\tquatretrois}{\begin{picture}(15,18)(-5,-1)
\put(3,0){\circle*{2}}
\put(-0.65,0){$\vee$}
\put(6,7){\circle*{2}}
\put(0,7){\circle*{2}}
\put(6,14){\circle*{2}}
\put(6,7){\line(0,1){7}}
\end{picture}}
\newcommand{\tquatrequatre}{\begin{picture}(15,18)(-5,-1)
\put(3,5){\circle*{2}}
\put(-0.65,5){$\vee$}
\put(6,12){\circle*{2}}
\put(0,12){\circle*{2}}
\put(3,0){\circle*{2}}
\put(3,0){\line(0,1){5}}
\end{picture}}
\newcommand{\tquatrecinq}{\begin{picture}(9,19)(-2,-1)
\put(0,0){\circle*{2}}
\put(0,0){\line(0,1){5}}
\put(0,5){\circle*{2}}
\put(0,5){\line(0,1){5}}
\put(0,10){\circle*{2}}
\put(0,10){\line(0,1){5}}
\put(0,15){\circle*{2}}
\end{picture}}

\newcommand{\tcinqun}{\begin{picture}(20,8)(-5,-1)
\put(3,0){\circle*{2}}
\put(-0.5,0){$\vee$}
\put(6,7){\circle*{2}}
\put(0,7){\circle*{2}}
\put(3,0){\line(2,1){10}}
\put(3,0){\line(-2,1){10}}
\put(-7,5){\circle*{2}}
\put(13,5){\circle*{2}}
\end{picture}}
\newcommand{\tcinqdeux}{\begin{picture}(15,14)(-5,-1)
\put(3,0){\circle*{2}}
\put(-0.65,0){$\vee$}
\put(6,7){\circle*{2}}
\put(0,7){\circle*{2}}
\put(3,7){\circle*{2}}
\put(3,0){\line(0,1){7}}
\put(0,7){\line(0,1){7}}
\put(0,14){\circle*{2}}
\end{picture}}
\newcommand{\tcinqtrois}{\begin{picture}(15,15)(-5,-1)
\put(3,0){\circle*{2}}
\put(-0.65,0){$\vee$}
\put(6,7){\circle*{2}}
\put(0,7){\circle*{2}}
\put(3,7){\circle*{2}}
\put(3,0){\line(0,1){7}}
\put(3,7){\line(0,1){7}}
\put(3,14){\circle*{2}}
\end{picture}}
\newcommand{\tcinqquatre}{\begin{picture}(15,14)(-5,-1)
\put(3,0){\circle*{2}}
\put(-0.65,0){$\vee$}
\put(6,7){\circle*{2}}
\put(0,7){\circle*{2}}
\put(3,7){\circle*{2}}
\put(3,0){\line(0,1){7}}
\put(6,7){\line(0,1){7}}
\put(6,14){\circle*{2}}
\end{picture}}
\newcommand{\tcinqcinq}{\begin{picture}(15,19)(-5,-1)
\put(3,0){\circle*{2}}
\put(-0.65,0){$\vee$}
\put(6,7){\circle*{2}}
\put(0,7){\circle*{2}}
\put(6,14){\circle*{2}}
\put(6,7){\line(0,1){7}}
\put(0,14){\circle*{2}}
\put(0,7){\line(0,1){7}}
\end{picture}}
\newcommand{\tcinqsix}{\begin{picture}(15,20)(-7,-1)
\put(3,0){\circle*{2}}
\put(-0.65,0){$\vee$}
\put(6,7){\circle*{2}}
\put(0,7){\circle*{2}}
\put(-3.65,7){$\vee$}
\put(3,14){\circle*{2}}
\put(-3,14){\circle*{2}}
\end{picture}}
\newcommand{\tcinqsept}{\begin{picture}(15,8)(-5,-1)
\put(3,0){\circle*{2}}
\put(-0.65,0){$\vee$}
\put(6,7){\circle*{2}}
\put(0,7){\circle*{2}}
\put(2.35,7){$\vee$}
\put(3,14){\circle*{2}}
\put(9,14){\circle*{2}}
\end{picture}}
\newcommand{\tcinqhuit}{\begin{picture}(15,26)(-5,-1)
\put(3,0){\circle*{2}}
\put(-0.65,0){$\vee$}
\put(6,7){\circle*{2}}
\put(0,7){\circle*{2}}
\put(0,14){\circle*{2}}
\put(0,7){\line(0,1){7}}
\put(0,21){\circle*{2}}
\put(0,14){\line(0,1){7}}
\end{picture}}
\newcommand{\tcinqneuf}{\begin{picture}(15,26)(-5,-1)
\put(3,0){\circle*{2}}
\put(-0.65,0){$\vee$}
\put(6,7){\circle*{2}}
\put(0,7){\circle*{2}}
\put(6,14){\circle*{2}}
\put(6,7){\line(0,1){7}}
\put(6,21){\circle*{2}}
\put(6,14){\line(0,1){7}}
\end{picture}}
\newcommand{\tcinqdix}{\begin{picture}(15,19)(-5,-1)
\put(3,5){\circle*{2}}
\put(-0.5,5){$\vee$}
\put(6,12){\circle*{2}}
\put(0,12){\circle*{2}}
\put(3,0){\circle*{2}}
\put(3,0){\line(0,1){12}}
\put(3,12){\circle*{2}}
\end{picture}}
\newcommand{\tcinqonze}{\begin{picture}(15,26)(-5,-1)
\put(3,5){\circle*{2}}
\put(-0.65,5){$\vee$}
\put(6,12){\circle*{2}}
\put(0,12){\circle*{2}}
\put(3,0){\circle*{2}}
\put(3,0){\line(0,1){5}}
\put(0,12){\line(0,1){7}}
\put(0,19){\circle*{2}}
\end{picture}}
\newcommand{\tcinqdouze}{\begin{picture}(15,26)(-5,-1)
\put(3,5){\circle*{2}}
\put(-0.65,5){$\vee$}
\put(6,12){\circle*{2}}
\put(0,12){\circle*{2}}
\put(3,0){\circle*{2}}
\put(3,0){\line(0,1){5}}
\put(6,12){\line(0,1){7}}
\put(6,19){\circle*{2}}
\end{picture}}
\newcommand{\tcinqtreize}{\begin{picture}(5,26)(-2,-1)
\put(0,0){\circle*{2}}
\put(0,0){\line(0,1){7}}
\put(0,7){\circle*{2}}
\put(0,7){\line(0,1){7}}
\put(0,14){\circle*{2}}
\put(-3.65,14){$\vee$}
\put(-3,21){\circle*{2}}
\put(3,21){\circle*{2}}
\end{picture}}
\newcommand{\tcinqquatorze}{\begin{picture}(9,26)(-5,-1)
\put(0,0){\circle*{2}}
\put(0,0){\line(0,1){5}}
\put(0,5){\circle*{2}}
\put(0,5){\line(0,1){5}}
\put(0,10){\circle*{2}}
\put(0,10){\line(0,1){5}}
\put(0,15){\circle*{2}}
\put(0,15){\line(0,1){5}}
\put(0,20){\circle*{2}}
\end{picture}}

\newcommand{\tdquatredeux}[4]{\begin{picture}(20,20)(-5,-1)
\put(3,0){\circle*{2}}
\put(-.65,0){$\vee$}
\put(6,7){\circle*{2}}
\put(0,7){\circle*{2}}
\put(0,14){\circle*{2}}
\put(0,7){\line(0,1){7}}
\put(5,-2){\tiny #1}
\put(9,5){\tiny #2}
\put(-5,5){\tiny #3}
\put(-5,12){\tiny #4}
\end{picture}}

\newcommand{\bun}{
\begin{picture}(5,5)(-2,0)
\put(0,0){\line(0,0){5}}
\end{picture}}

\newcommand{\bdeux}{
\begin{picture}(10,15)(-5,0)
\put(-3.5,4.1){$\vee$}
\put(0,0){\line(0,1){5}}
\end{picture}}

\newcommand{\btroisun}{
\begin{picture}(15,20)(-7,0)
\put(-3.5,4.1){$\vee$}
\put(0,0){\line(0,1){5}}
\put(-6.3,10.7){$\vee$}
\end{picture}}
\newcommand{\btroisdeux}{
\begin{picture}(15,20)(-5,0)
\put(-3.5,4.1){$\vee$}
\put(0,0){\line(0,1){5}}
\put(-0.7,10.7){$\vee$}
\end{picture}}

\newcommand{\bquatreun}{
\begin{picture}(15,25)(-7,0)
\put(-3.5,4.1){$\vee$}
\put(0,0){\line(0,1){5}}
\put(-6.3,10.7){$\vee$}
\put(-9.1,17.3){$\vee$}
\end{picture}}
\newcommand{\bquatredeux}{
\begin{picture}(15,25)(-7,0)
\put(-3.5,4.1){$\vee$}
\put(0,0){\line(0,1){5}}
\put(-6.3,10.7){$\vee$}
\put(-3.5,17.3){$\vee$}
\end{picture}}
\newcommand{\bquatretrois}{
\begin{picture}(15,25)(-7,0)
\put(-3.5,4.1){$\vee$}
\put(0,0){\line(0,1){5}}
\put(-0.7,10.7){$\vee$}
\put(-3.5,17.3){$\vee$}
\end{picture}}
\newcommand{\bquatrequatre}{
\begin{picture}(15,25)(-7,0)
\put(-3.5,4.1){$\vee$}
\put(0,0){\line(0,1){5}}
\put(-0.7,10.7){$\vee$}
\put(2.1,17.3){$\vee$}
\end{picture}}
\newcommand{\bquatrecinq}{
\begin{picture}(15,20)(-7,0)
\put(0,0){\line(0,1){5}}
\put(0,5){\line(2,1){10}}
\put(0,5){\line(-2,1){10}}
\put(-3,6.5){\line(0,1){5}}
\put(3,6.5){\line(0,1){5}}
\end{picture}}

\input{xy}
\xyoption{all}

\newcommand{\tdelta}{\tilde{\Delta}}
\newcommand{\h}{{\cal H}}
\newcommand{\F}{\mathbf{F}}
\newcommand{\T}{\mathbf{T}}

\title{The infinitesimal Hopf algebra and the poset of planar forests}
\date{}
\author{L. Foissy \\
\\
{\small{\it Laboratoire de Math\'ematiques, Universit\'e de Reims}}\\
\small{{\it Moulin de la Housse - BP 1039 - 51687 REIMS Cedex 2, France}}\\
\small{e-mail : loic.foissy@univ-reims.fr}}

\newtheorem{defi}{\indent Definition}
\newtheorem{lemma}[defi]{\indent Lemma}
\newtheorem{cor}[defi]{\indent Corollary}
\newtheorem{theo}[defi]{\indent Theorem}
\newtheorem{prop}[defi]{\indent Proposition}

\begin{document}

\maketitle

ABSTRACT We introduce an infinitesimal Hopf algebra of planar trees, 
generalising the construction of the non-commutative Connes-Kreimer Hopf algebra.
A non-degenerate pairing and a dual basis are defined, and a combinatorial interpretation
of the pairing in terms of orders on the vertices of planar forests is given.
Moreover, the coproduct and the pairing can also be described with the help of a partial
order on the set of planar forests, making it isomorphic to the Tamari poset.
As a corollary, the dual basis can be computed with a M\"obius inversion.\\

KEY-WORDS Infinitesimal Hopf algebra, planar tree, Tamari poset.\\

AMS CLASSIFICATION 16W30, 05C05, 06A11

\tableofcontents

\section*{Introduction}

The Connes-Kreimer Hopf algebra of rooted trees is introduced and studied in \cite{Connes,Foissy,Kreimer1,Kreimer2,Kreimer3}.
This commutative, non commutative Hopf algebra is used to treat a problem of Renormalisation in Quantum Fields Theory,
as explained in \cite{Connes2,Connes3}. A non-commutative version of this Hopf algebra is introduced simultaneously in \cite{Foissy2} and \cite{Holtkamp}.
This Hopf algebra $\h_{P,R}$, based on planar rooted trees, is neither commutative nor cocommutative, and satisfies a universal property
in Cartier-Quillen cohomology. This property is used in \cite{Foissy} to prove that $\h_{P,R}$ is isomorphic to its (graded) dual.
In other terms, $\h_{P,R}$ owns a non-degenerate, symmetric Hopf pairing, and a dual basis of its basis of planar forests.
This pairing admits a description in terms of two partial orders on the vertices of the planar forests.

Our aim in the present text is to introduce an infinitesimal version of this Hopf algebra $\h_{P,R}$. 
The concept of infinitesimal Hopf algebra is introduced in \cite{Loday2}. Namely, an infinitesimal  bialgebra is a space $A$,
both an associative, unitary algebra and a coassociative, counitary coalgebra, with the following compatibility:
$$\Delta(ab)=(a\otimes 1)\Delta(b)+\Delta(a)(1\otimes b)-a\otimes b.$$
If it has an antipode, $A$ will be said to be an infinitesimal Hopf algebra.
It is proved in \cite{Loday2} that an infinitesimal bialgebra $A$, which is connected as a coalgebra, is isomorphic
to $T(Prim(A))$, with its concatenation product and deconcatenation coproduct: this is the infinitesimal rigidity theorem.

We here construct an infinitesimal coproduct over the algebra $\h$ of planar rooted trees (theorem \ref{9}).
We use for this the fact that $\h$, given the linear application $B^+$ of grafting on a common root, is an initial object in a certain category.
This infinitesimal coproduct is given by left-admissible cuts (theorem \ref{10}), whereas the usual Hopf coproduct is given by admissible cuts.
We also give a description of this coproduct in terms of the two partial orders $\geq_{high}$ and $\geq_{left}$ on the vertices of a planar forest
(proposition \ref{12}). We also give a formula for the antipode in terms of left cuts (proposition \ref{15}). Using the infinitesimal universal property of $\h$ 
(theorem \ref{16}), we construct a non-degenerate Hopf pairing between $\h$ and $\h^{op,cop}$ (theorem \ref{19}), 
and a dual basis $(f_F)$ of the basis of forests of $\h$. This pairing $\langle-,-\rangle$ admits a combinatorial interpretation, described in theorem \ref{24}.
All these results are infinitesimal versions of the classical Hopf results of \cite{Foissy2}.

Differences between the infinitesimal and the Hopf case become clear with the observation that the pairing
$\langle F,G\rangle$ of two forests $F$ and $G$ is always $0$ or $1$ in the infinitesimal case. 
This leads to an interpretation of this pairing in terms of a certain poset, namely the poset of planar forest.
A partial order is defined on the set of planar forests with the help of certain transformations of forests (definition \ref{25}).
This poset $\F$ is isomorphic to the Tamari poset of planar binary trees \cite{Stanley2}, as it is proved in theorem \ref{31}. 
As a consequence, it has a decreasing isomorphism $m$, corresponding to the vertical symmetry of planar binary trees in the Tamari poset.
The pairing $\langle-,-\rangle$ satisfies the following property: for all planar forests $F$ and $G$, $\langle F,G\rangle=1$ if, and only if,
$F\geq m(G)$  in $\F$. As a consequence, the dual basis $(f_F)$ is given by a M\"obius inversion (corollary \ref{39}).
Moreover, the product of two elements of the dual basis admits also a description using suborders of $\F$ (corollary \ref{29}).
We shall show in another text that this dual basis can be iteratively computed with the help of two operads of planar forests.
For the sake of simplicity, we restrict ourselves here to planar rooted trees with no decorations, 
but  there exists versions of all these results for planar decorated rooted trees, and these versions are proved similarly.\\

This paper is organised as follows: the first paragraph is devoted to recalls and complements about infinitesimal Hopf algebras.
In particular, it is proved that for any infinitesimal Hopf algebra $A$, $Ker(\varepsilon)=Prim(A)\oplus Ker(\varepsilon)^2$, and the projector
on $Prim(A)$ in this direct sum is the antipode, recovering in this way the rigidity theorem of \cite{Loday2}.
The infinitesimal Hopf algebra of planar rooted trees $\h$ is introduced in the second section. We construct its infinitesimal coproduct
and give its description in terms of left-admissible cuts and biideals, before a formula for the antipode.
We prove a universal property of $\h$ and use it to construct a Hopf pairing between $\h$ and $\h^{op,cop}$.
The combinatorial description of this pairing is then given. The last section deals with the poset of forests $\F$ and its applications.
We prove that this poset is isomorphic to the Tamari poset and describe a decreasing isomorphism of $\F$. 
The link between the pairing of $\h$ and the order on $\F$ is then given.\\

{\bf Notation.} We denote by $K$ a commutative field, of any characteristic. Every algebra, coalgebra, etc, will be taken over $K$.
\section{Recalls on infinitesimal Hopf algebras}

We refer to \cite{Abe,Sweedler} for the classical results and definitions about coalgebras, bialgebras, Hopf algebras.

\subsection{Infinitesimal Hopf algebras}

\begin{defi}
\textnormal{(See \cite{Loday2}).
\begin{enumerate}
\item An {\it infinitesimal bialgebra} is an associative, unitary algebra $A$, together with a coassociative, counitary coproduct, 
satisfying the following compatibility: for all $a,b \in A$,
\begin{equation}
\label{E1} \Delta(ab)=\Delta(a)(1\otimes b)+(a\otimes 1)\Delta(b)-a\otimes b.
\end{equation}
\item Let $A$ be an infinitesimal bialgebra. If $Id_A$ has an inverse $S$ in the associative convolution algebra $({\cal L}(A),\star)$, 
we shall say that $A$ is an {\it infinitesimal Hopf algebra}, and $S$ will be called the {\it antipode} of $A$.
\end{enumerate}}
\end{defi}

{\bf Remarks.}
\begin{enumerate}
\item This is not the same definition as used by Aguiar in \cite{Aguiar3}.
\item Let $A$ be an infinitesimal bialgebra and let $M$ be the kernel of its counit. We shall prove in proposition \ref{2} that  $M$ is an ideal.
Moreover, $M$ is given a coassociative, non counitary coproduct $\tdelta$ defined by:
$$\tdelta: \left\{ \begin{array}{rcl}
M&\longrightarrow & M\otimes M\\
x&\longrightarrow &\Delta(x)-x\otimes 1-1\otimes x.
\end{array}
\right.$$
The compatibility between $\tdelta$ and the product is given by the {\it non (co)unital infinitesimal compatibility}:
$$\tdelta(ab)=(a\otimes 1)\tdelta(b)+\tdelta(a)(1\otimes b)+a\otimes b.$$
\item By induction, if $x_1,\ldots, x_n\in M$:
$$\tdelta(x_1\ldots x_n)=\sum_{i=1}^{n-1} x_1 \ldots x_i \otimes x_{i+1}\ldots x_n
+\sum_{i=1}^n (x_1 \ldots x_{i-1} \otimes 1)\tdelta(x_i) (1\otimes x_{i+1}\ldots x_n).$$
\item In particular, if $x_1,\ldots,x_n$ are primitive elements of $A$:
$$\Delta(x_1\ldots x_n)=\sum_{i=0}^n x_1\ldots x_i \otimes x_{i+1}\ldots x_n.$$
\end{enumerate}

{\bf Notations.} Let $A$ be an infinitesimal bialgebra. For all $x\in A$, we denote $\Delta(x)=x^{(1)}\otimes x^{(2)}$. 
Moreover, if $\varepsilon(x)=0$, we denote $\tdelta(x)=x'\otimes x''$.\\

{\bf Examples.}
\begin{enumerate}
\item Let $V$ be a vector space. The tensor algebra $T(V)$ is given a structure of infinitesimal Hopf algebra with the coproduct $\Delta$ defined,
for $v_1,\ldots,v_n \in V$, by:
$$\Delta(v_1\ldots v_n)=\sum_{i=0}^n v_1\ldots v_i \otimes v_{i+1}\ldots v_n.$$
It is proved in \cite{Loday2} that any connected (as a coalgebra) infinitesimal bialgebra $A$ is isomorphic to $T(Prim(A))$.
\item If $A$ is an infinitesimal Hopf algebra, then $A^{op,cop}$ also is, with the same antipode. 
Note that $A^{op}$ and $A^{cop}$ are not infinitesimal bialgebras, as the compatibility (\ref{E1}) is no more satisfied.
\item If $A$ is a graded infinitesimal Hopf algebra, such that its homogeneous components are finite-dimensional, then its graded dual $A^*$ also is.
\end{enumerate}

\begin{prop}
\label{2}
Let $A$ be an infinitesimal bialgebra.
\begin{enumerate}
\item Then $\Delta(1)=1\otimes 1$. In other terms, the unit application $\nu$ is a coalgebra morphism:
$$\nu: \left\{ \begin{array}{rcl}
K&\longrightarrow &A \\
1&\longrightarrow &1.
\end{array}
\right.$$
\item The counit $\varepsilon$ is an algebra morphism.
\item If, moreover, $A$ is an infinitesimal Hopf algebra, then $S(1)=1$ and $\varepsilon\circ S=\varepsilon$.
\end{enumerate}
\end{prop}

{\bf Proof.}\begin{enumerate}
\item For $a=b=1$, relation (\ref{E1}) becomes $\Delta(1)=\Delta(1)+\Delta(1)-1\otimes 1$. So $\Delta(1)=1\otimes 1$.
As a consequence, $\varepsilon(1)=1$. Moreover, if $A$ has an antipode, $S(1)=1$.
\item For $a,b \in A$:
\begin{eqnarray*}
(\varepsilon \otimes \varepsilon)\circ \Delta(ab)&=&\varepsilon(a^{(1)})\varepsilon(a^{(2)}b)
+\varepsilon(ab^{(1)}) \varepsilon(b^{(2)})-\varepsilon(a)\varepsilon(b)\\
&=&\varepsilon(ab)+\varepsilon(ab)-\varepsilon(a)\varepsilon(b)\\
&=&\varepsilon(ab).
\end{eqnarray*}
So $\varepsilon(ab)=\varepsilon(a)\varepsilon(b)$. 
\item For all $a\in A$:
$$\varepsilon(a)=\varepsilon(\varepsilon(a)1)=\varepsilon\circ m\circ (S \otimes Id) \circ \Delta(a)
=\varepsilon\left(S\left(a^{(1)}\right)\right) \varepsilon\left(a^{(2)}\right)=\varepsilon(S(a)).$$
So $\varepsilon \circ S=\varepsilon$. $\Box$
\end{enumerate}

\begin{lemma}
\label{3}
\begin{enumerate}
\item Let $A$, $B$ be two augmented algebras, with respective augmentations denoted by 
$\varepsilon_A:A\longrightarrow K$ and $\varepsilon_B:B\longrightarrow K$. 
Then $A\otimes B$ is an associative, unitary algebra,  with product $._{\varepsilon_{A,B}}$ 
given for all $a_1,a_2\in A$, $b_1,b_2\in B$, by:
$$(a_1\otimes b_1)._{\varepsilon_{A,B}} (a_2 \otimes b_2)=
\varepsilon(a_2) a_1\otimes b_1 b_2+\varepsilon(b_1) a_1a_2\otimes b_2-\varepsilon_A(a_2)\varepsilon_B(b_1)a_1\otimes b_2.$$
The unit is $1_A\otimes 1_B$.
\item Let $A$, $B$ be two pointed coalgebras, with group-like elements $1_A$ and $1_B$. Then $A\otimes B$ is a coassociative, counitary coalgebra, 
with coproduct $\Delta_{1_{A,B}}$ given by:
$$\Delta_{1_{A,B}}(a \otimes b)=a\otimes b^{(1)} \otimes 1_A \otimes b^{(2)}
+a^{(1)} \otimes 1_B \otimes a^{(2)} \otimes b-a \otimes 1_B \otimes 1_A \otimes b.$$
The counit is $\varepsilon_A \otimes \varepsilon_B$.
\end{enumerate}
\end{lemma}

{\bf Proof.} Direct computations. $\Box$\\

{\bf Remarks.} \begin{enumerate}
\item When the augmented algebras $A$ and $B$ are equal, we shall simply denote $\varepsilon_A=\varepsilon_B=\varepsilon$, 
and $\Delta_{\varepsilon_{A,B}}=\Delta_\varepsilon$. 
When the pointed coalgebras $A$ and $B$ are equal, we shall denote $1_A=1_B=1$ and $\Delta_{1_{A,B}}=\Delta_1$. 
\item If $A$ is an infinitesimal bialgebra, then compatibility (\ref{E1}) means that the coproduct $\Delta:(A,.) \longrightarrow (A\otimes A,._\varepsilon)$,
where $\varepsilon$ is the counit of $A$, is a morphism of algebras. Indeed, for all $a,b \in A$:
\begin{eqnarray*}
\Delta(a)._\varepsilon \Delta(b)&=&\varepsilon(b^{(1)}) a^{(1)} \otimes a^{(2)}b^{(2)}+\varepsilon(a^{(2)})a^{(1)}b^{(1)} \otimes b^{(2)}
-\varepsilon(a^{(2)})\varepsilon(b^{(1)})a^{(1)} \otimes b^{(2)}\\
&=&a^{(1)} \otimes a^{(2)}b+ab^{(1)} \otimes b^{(2)}-a\otimes b.
\end{eqnarray*}

Dually, it also means that the product $m:(A\otimes A,\Delta_1) \longrightarrow (A,\Delta)$
is a morphism of coalgebras.
\end{enumerate}

\subsection{Antipode of an infinitesimal Hopf algebra}

\begin{lemma}
\label{4}
Let $A$ be an infinitesimal Hopf algebra. 
\begin{enumerate}
\item For all $a,b \in A$, $S(ab)=\varepsilon(a) S(b)+\varepsilon(b) S(a)-\varepsilon(a)\varepsilon(b)1$.
 In particular, for all $a,b \in A$, such that  $\varepsilon(a)=\varepsilon(b)=0$, $S(ab)=0$.
\item For all $a\in A$, $\Delta(S(a))=S(a)\otimes 1+1\otimes S(a)-\varepsilon(a) 1\otimes 1$.
In particular, for all $a\in A$, such that $\varepsilon(a)=0$, $S(a)$ is primitive.
\end{enumerate}
\end{lemma}

{\bf Proof.} \begin{enumerate}
\item Let us consider the convolution algebra ${\cal L}(A\otimes A,A)$, where $A\otimes A$ is given the coproduct of lemma \ref{3}. For all $a,b \in A$:
\begin{eqnarray*}
&&((S \circ m) \star m)(a \otimes b)\\
&=&m\circ ((S\circ m)\otimes m)\circ \Delta_1(a\otimes b)\\
&=&m\circ ((S\circ m)\otimes m)\left(a^{(1)} \otimes 1\otimes a^{(2)}\otimes b
+a\otimes b^{(1)} \otimes 1 \otimes b^{(2)}-a\otimes 1\otimes 1\otimes b\right)\\ 
&=&S(a^{(1)})a^{(2)}b+S(ab^{(1)})b^{(2)}-S(a)b\\
&=&m\circ (S \otimes Id)\circ \Delta(ab)\\
&=&\varepsilon(ab)1\\
&=&\varepsilon(a)\varepsilon(b)1.
\end{eqnarray*}
So $S \circ m$ is a left inverse of $m$.

Let $T:A\otimes A\longrightarrow A$ defined by $T(a\otimes b)=\varepsilon(a)S(b)+\varepsilon(b)S(a)-\varepsilon(a)\varepsilon(b)1$.
Let us compute $m\star T$ in ${\cal L}(A\otimes A,A)$:
\begin{eqnarray*}
&&(m\star T)(a\otimes b)\\
&=&m\circ (m\otimes T)\circ \Delta_1(a\otimes b)\\
&=&m\circ (m\otimes T)\left(a^{(1)} \otimes 1\otimes a^{(2)}\otimes b
+a\otimes b^{(1)} \otimes 1 \otimes b^{(2)}-a\otimes 1\otimes 1\otimes b\right)\\ 
&=&a^{(1)}\left(\varepsilon(a^{(2)}) S(b)+\varepsilon(b) S(a^{(2)})-\varepsilon(a^{(2)}) \varepsilon(b)1 \right)\\
&&+a b^{(1)}\left(S(b^{(2)})+\varepsilon(b^{(2)}) 1-\varepsilon(b^{(2)})1 \right)-a\left(S(b)+\varepsilon(b)1-\varepsilon(b)1 \right)\\
&=&aS(b)+\varepsilon(a)\varepsilon(b)1-\varepsilon(b)a+\varepsilon(b)a-aS(b)\\
&=&\varepsilon(a)\varepsilon(b)1.
\end{eqnarray*}
So $T$ is a right inverse of $m$. As the convolution product is associative, $S \circ m=T$.
\item Let us consider the convolution algebra ${\cal L}(A,A\otimes A)$, where $A\otimes A$ is given the product of lemma \ref{3}. For all $a \in A$:
\begin{eqnarray*}
&&((\Delta \circ S) \star \Delta)(a)\\
&=& m_\varepsilon \circ ((\Delta \circ S) \otimes \Delta) \circ \Delta (a)\\
&=&m_\varepsilon \left(S(a^{(1)})^{(1)}\otimes S(a^{(1)})^{(2)} \otimes a^{(2)} \otimes a^{(3)} \right)\\
&=&\varepsilon(S(a^{(1)})^{(2)}) S(a^{(1)})^{(1)} a^{(2)} \otimes a^{(3)}+\varepsilon(a^{(2)})S(a^{(1)})^{(1)}\otimes S(a^{(1)})^{(2)}a^{(3)}\\
&&-\varepsilon(S(a^{(1)})^{(2)}) \varepsilon(a^{(2)})S(a^{(1)})^{(1)}\otimes a^{(3)}\\
&=&S(a^{(1)})a^{(2)} \otimes a^{(3)}+S(a^{(1)})^{(1)}\otimes S(a^{(1)})^{(2)}a^{(2)}-S(a^{(1)})\otimes a^{(2)}\\
&=&\Delta(S(a^{(1)})a^{(2)})\\
&=&\varepsilon(a)\Delta(1)\\
&=&\varepsilon(a)1\otimes 1.
\end{eqnarray*}
So $\Delta \circ S$ is a left inverse of $\Delta$.

Let $T:A\longrightarrow A\otimes A$ defined by $T(a)=S(a)\otimes 1+1\otimes S(a)-\varepsilon(a)1\otimes 1$.
Let us compute $\Delta \star T$ in ${\cal L}(A,A\otimes A)$:
\begin{eqnarray*}
&&(\Delta \star T)(a)\\
&=&m_\varepsilon \circ (\Delta \otimes T) \circ \Delta(a)\\
&=&m_\varepsilon (a^{(1)} \otimes a^{(2)}\otimes 1\otimes S(a^{(3)})+a^{(1)}\otimes a^{(2)}\otimes S(a^{(3)})\otimes 1
-\varepsilon(a^{(3)})a^{(1)}\otimes a^{(2)} \otimes 1\otimes 1)\\
&=&a^{(1)} \otimes a^{(2)}S(a^{(3)})+\varepsilon(a^{(2)})a^{(1)}\otimes S(a^{(3)})-\varepsilon(a^{(2)})a^{(1)}\otimes S(a^{(3)})\\
&&+\varepsilon(S(a^{(3)}))a^{(1)}\otimes a^{(2)}+\varepsilon(a^{(2)})a^{(1)}S(a^{(3)})\otimes 1-\varepsilon(a^{(2)})\varepsilon(S(a^{(3)}))a^{(1)}\otimes 1\\
&&-a^{(1)}\otimes a^{(2)}-\varepsilon(a^{(2)})a^{(1)}\otimes 1+\varepsilon(a^{(2)})a^{(1)}\otimes 1\\
&=&a\otimes 1+a^{(1)}\otimes a^{(2)}+\varepsilon(a)1\otimes 1-a\otimes 1-a^{(1)}\otimes a^{(2)}-a\otimes 1+a\otimes 1\\
&=&\varepsilon(a)1\otimes 1.
\end{eqnarray*}
So $T$ is a right inverse of $\Delta$. As the convolution product is associative, $\Delta \circ S=T$. $\Box$
\end{enumerate}

\begin{cor}
Let $A$ be an infinitesimal Hopf algebra. Then $Ker(\varepsilon)=Prim(A)\oplus Ker(\varepsilon)^2$.
The projection on $Prim(A)$ in this direct sum is $-S$.
\end{cor}

{\bf Proof.} Let $a\in Ker(\varepsilon)$. Then $\Delta(a)=a\otimes 1+1\otimes a +a'\otimes a''$, 
with $a'\otimes a'' \in Ker(\varepsilon) \otimes Ker(\varepsilon)$. Moreover:
$$0=\varepsilon(a) 1=m\circ (S \otimes Id) \circ \Delta(a)=S(a)+a+S(a')a'',$$
so $a=-S(a)-S(a')a''$. By lemma \ref{4}, $-S(a) \in Prim(A)$ and $S(a')a'' \in S(Ker(\varepsilon))Ker(\varepsilon)^2\subseteq Ker(\varepsilon)^2$,
so $Ker(\varepsilon)=Prim(A)+Ker(\varepsilon)^2$. If $a\in Prim(A)$, then $a'\otimes a''=0$, so $-S(a)=a$. Moreover, $S(Ker(\varepsilon)^2)=(0)$,
so $Ker(\varepsilon)=Prim(A)\oplus Ker(\varepsilon)^2$ and the projector on $Prim(A)$ in this direct sum is $-S$. $\Box$\\

{\bf Remarks.} \begin{enumerate}
\item This result implies the rigidity theorem of \cite{Loday2}.
\item It is also possible to prove lemma \ref{4} using braided Hopf algebras.
\end{enumerate}

\section{Infinitesimal Hopf algebra of planar trees}

\subsection{Algebra of planar trees and universal property}

\begin{defi}
\textnormal{\begin{enumerate}
\item The set of planar rooted trees will be denoted by $\T$ (see \cite{Foissy2,Holtkamp}). 
\item The algebra $\h$ is the free associative algebra generated by $\T$. The monomials of $\h$ will be called {\it planar forests}. 
The set of planar forests will be denoted by $\F$. The {\it weight} of an element $F \in \F$ is the number of its vertices. 
\end{enumerate}}
\end{defi}

{\bf Examples.} \begin{enumerate}
\item Planar rooted trees of weight $\leq 5$:
$$\tun,\tdeux,\ttroisun,\ttroisdeux,\tquatreun, \tquatredeux,\tquatretrois,\tquatrequatre,\tquatrecinq,\tcinqun,\tcinqdeux,\tcinqtrois,\tcinqquatre,\tcinqcinq,
\tcinqsix,\tcinqsept,\tcinqhuit,\tcinqneuf,\tcinqdix,\tcinqonze,\tcinqdouze,\tcinqtreize,\tcinqquatorze.$$
\item Planar rooted forests of weight $\leq 4$:
$$1,\tun,\tun\tun,\tdeux,\tun\tun\tun,\tdeux\tun,\tun \tdeux,\ttroisun,\ttroisdeux,\tun\tun\tun\tun,\tdeux\tun\tun,\tun \tdeux \tun, \tun \tun \tdeux,
\ttroisun\tun,\tun \ttroisun,\ttroisdeux\tun,\tun \ttroisdeux,\tdeux\tdeux,\tquatreun,\tquatredeux,\tquatretrois,\tquatrequatre,\tquatrecinq.$$
\end{enumerate}

We define the operator $B^+:\h\longrightarrow \h$, which associates, to a forest $F\in \F$, the tree obtained by grafting the roots of the trees of $F$ 
on a common root. For example, $B^+(\tdeux\tun)=\tquatredeux$, and $B^+(\tun\tdeux)=\tquatretrois$. It is shown in \cite{Moerdijk} 
that $(\h,B^+)$ is an initial object in the category of couples $(A,L)$, where $A$ is an algebra, and $L:A \longrightarrow A$ any linear operator.
More explicitely:

\begin{theo}[Universal property of $\h$]
Let $A$ be any algebra and let $L:A\longrightarrow A$ be a linear map. 
Then there exists a unique algebra morphism $\phi:\h \longrightarrow A$, such that $\phi \circ B^+=L\circ \phi$.
\end{theo}

{\bf Remark.} Note that $\phi$ is inductively defined in the following way: for all trees $t_1,\ldots, t_n \in \T$,
$$  \left\{ \begin{array}{rcl}
\phi(1)&=&1,\\
\phi(t_1\ldots t_n)&=&\phi(t_1)\ldots \phi(t_n),\\
\phi(B^+(t_1\ldots t_n))&=&L(\phi(t_1)\ldots \phi(t_n)).
\end{array}
\right.$$

The end of this paragraph is devoted to the introduction of several combinatorial concepts, which will be useful for the sequel.

\begin{defi}
\textnormal{Let $F \in \F$. An {\it admissible cut} is a non empty cut of certain edges and trees of $F$, 
such that each path in a non-cut tree of $F$ meets at most one cut edge (see \cite{Connes,Foissy2}).
The set of admissible cuts of $F$ will be denoted by $Adm(F)$.
If $c$ is an admissible cut of $F$, the forest of the vertices which are over the cuts of $c$ will be denoted
by $P^c(t)$ (branch of the cut $c$), and the remaining forest will be denoted by $R^c(t)$ (trunk of the cut).}
\end{defi}

We now recall several order relations on the set of the vertices of a planar forest, see \cite{Foissy2} for more details.
Let $F=t_1\ldots t_n \in \F-\{1\}$ and let $s,s'$ be two vertices of $F$.
\begin{enumerate}
\item  We shall say that $s \geq_{high} s'$ if there exists a path from $s'$ to $s$ in $F$, the edges of $F$ being oriented from the roots to the leaves. 
Note that $\geq_{high}$ is a partial order, whose Hasse graph is the forest $F$.
\item If $s$ and $s'$ are not comparable for $\geq_{high}$, we shall say that $s \geq_{left} s'$ if one of these assertions is satisfied:
\begin{enumerate}
\item $s$ is a vertex of $t_i$ and $s'$ is a vertex of $t_j$, with $i<j$.
\item $s$ and $s'$ are vertices of the same $t_i$, and $s \geq_{left} s'$ in the forest obtained from $t_i$ by deleting its root.
\end{enumerate}
This defines the partial order $\geq_{left}$ for all forests $F$, by induction on the the weight.
\item We shall say that $s \geq_{h,l} s'$ if $s \geq_{high} s'$ or $s \geq_{left} s'$. This defines a total order on the vertices of $F$. 
\end{enumerate}

{\bf Example.} Let $t=\tquatredeux$. We index its vertices in the following way: $\tdquatredeux{4}{3}{2}{1}$.
The following arrays give the order relations $\geq_{high}$ and $\geq_{left}$ for the vertices of $t$. 
A symbol $\times$ means that the vertices are not comparable for the order.

$$\begin{array}{|c|c|c|c|c|}
\hline x\setminus y&s_1&s_2&s_3&s_4\\
\hline s_1&=&\geq_{high}&\times&\geq_{high}\\
\hline s_2&\leq_{high}&=&\times&\geq_{high}\\
\hline s_3&\times&\times&=&\geq_{high}\\
\hline s_4&\leq_{high}&\leq_{high}&\leq_{high}&=\\
\hline 
\end{array} \hspace{1cm}
\begin{array}{|c|c|c|c|c|}
\hline x\setminus y&s_1&s_2&s_3&s_4\\
\hline s_1&=&\times&\geq_{left}&\times\\
\hline s_2&\times&=&\geq_{left}&\times\\
\hline s_3&\leq_{left}&\leq_{left}&=&\times\\
\hline s_4&\times&\times&\times&=\\
\hline 
\end{array} $$
So $s_1 \geq_{h,l} s_2 \geq_{h,l} s_3 \geq_{h,l} s_4$.

\subsection{Infinitesimal coproduct of $\h$}

We define:
$$ \varepsilon: \left\{ \begin{array}{rcl}
\h &\longrightarrow & K\\
F \in \F&\longrightarrow & \delta_{F,1}.
\end{array}\right.$$
Then $\varepsilon$ is clearly an algebra morphism. Moreover, $\varepsilon \circ B^+=0$. We also consider:
$$\nu: \left\{ \begin{array}{rcl}
K&\longrightarrow &\h \\
\lambda &\longrightarrow & \lambda 1.
\end{array}
\right.$$
Note that $\varepsilon \circ \nu=Id_K$.

\begin{theo}
\label{9}
Let $\Delta:\h \longrightarrow (\h \otimes \h,._\varepsilon)$ be the unique algebra morphism such that:
$$\Delta \circ B^+=(Id \otimes B^+ + B^+ \otimes (\nu \circ\varepsilon))\circ \Delta.$$
Then $(\h,\Delta)$ is an infinitesimal bialgebra. It is graded by the weight.
\end{theo}

{\bf Proof.} 

{\it First step.} Let us show that $\varepsilon$ is a counit for $\Delta$. Let $\varphi=(Id \otimes \varepsilon) \circ  \Delta$. 
By composition, $\varphi$ is an algebra endomorphism of $\h$. Moreover:
\begin{eqnarray*}
\varphi \circ B^+&=&(Id \otimes \varepsilon) \circ (Id \otimes B^+ + B^+ \otimes (\nu \circ\varepsilon))\circ \Delta\\
&=&(Id \otimes (\varepsilon \circ B^+)+B^+ \otimes (\varepsilon \circ \nu \circ \varepsilon))\circ \Delta\\
&=&B^+\circ (Id \otimes \varepsilon) \circ \Delta\\
&=&B^+ \circ \varphi. 
\end{eqnarray*}
By unicity in the universal property, $\varphi=Id_\h$. So $\varepsilon$ is a right counity for $\Delta$.
Let $\phi=(\varepsilon \otimes Id) \circ \Delta$. Then:
\begin{eqnarray*}
\phi \circ B^+&=&(\varepsilon \otimes Id) \circ (Id \otimes B^+ + B^+ \otimes (\nu \circ\varepsilon))\circ \Delta\\
&=&(\varepsilon \otimes B^+ +(\varepsilon \circ B^+) \otimes (\nu \circ\varepsilon))\circ \Delta\\
&=&B^+\circ (\varepsilon \otimes Id) \circ \Delta+0\\
&=&B^+ \circ \phi. 
\end{eqnarray*}
By unicity in the universal property, $\phi=Id_\h$. So $\varepsilon$ is a counit for $\Delta$.
As a consequence, $(\varepsilon \otimes \varepsilon) \circ \Delta=\varepsilon$.\\ 

{\it Second step.} Let us show that $\Delta$ is coassociative. We consider $\theta=(\Delta \otimes Id)\circ \Delta$. 
This is an algebra morphism from $\h$ to $\h \otimes \h \otimes \h$. Moreover:
\begin{eqnarray*}
\theta \circ B^+&=&
\left(\Delta	\otimes B^+ + (\Delta \circ B^+) \otimes (\nu \circ\varepsilon)\right)\circ \Delta\\
&=&\left(Id \otimes Id \otimes B^+ +Id \otimes B^+ \otimes (\nu \circ \varepsilon)
+B^+ \otimes (\nu \circ \varepsilon) \otimes (\nu \circ \varepsilon)\right)\circ \theta.
\end{eqnarray*}
Consider now $\theta'=(Id \otimes \Delta) \circ \Delta$. This is also an algebra morphism from $\h$ to $\h \otimes \h \otimes \h$. Moreover:
\begin{eqnarray*}
\theta' \circ B^+&=&
\left(Id \otimes (\Delta \circ B^+)
+B^+ \otimes (\Delta \circ \nu \circ \varepsilon)\right) \circ \Delta\\ 
&=& \left( Id \otimes Id \otimes B^++Id \otimes B^+ \otimes (\nu \circ \varepsilon) \right) \circ \theta'
+\left(B^+ \otimes ((\nu \otimes \nu)\circ \varepsilon) \right)\circ \Delta\\
&=& \left( Id \otimes Id \otimes B^++Id \otimes B^+ \otimes (\nu \circ \varepsilon) \right) \circ \theta'
+\left(B^+ \otimes (\nu \circ \varepsilon) \otimes (\nu\circ \varepsilon) \right)\circ \theta'\\
&=& \left( Id \otimes Id \otimes B^++Id \otimes B^+ \otimes (\nu \circ \varepsilon) +
B^+ \otimes (\nu \circ \varepsilon) \otimes (\nu\circ \varepsilon) \right)\circ \theta'.
\end{eqnarray*}
By unicity in the universal property, $\theta=\theta'$, so $\Delta$ is coassociative.
As $\Delta:\h \longrightarrow (\h\otimes \h,._\varepsilon)$ is a morphism of algebras, $(\h, \Delta)$ is an infinitesimal bialgebra.\\

{\it Last step.} It remains to show that $\Delta$ is homogeneous of degree $0$. 
Easy induction, using the fact that $Id \otimes B^+ + B^+ \otimes (\nu \circ\varepsilon)$ is homogeneous of degree $1$. 
Note that it can also be proved from proposition \ref{10}. $\Box$\\

{\bf Remarks.} \begin{enumerate}
\item In other terms, the coproduct $\Delta$ is uniquely defined by the following relations: for all $x,y \in \h$,
$$ \left\{ \begin{array}{rcl}
\Delta(1)&=&1\otimes 1,\\
\Delta(xy)&=&(x\otimes 1)\Delta(y)+\Delta(x)(1\otimes y)-x \otimes y,\\
\Delta(B^+(x))&=&(Id \otimes B^+) \circ \Delta(x) +B^+(x) \otimes 1.
\end{array} \right.$$
\item Equivalently, the non unitary coproduct $\tdelta$ satisfies the following property: for all $x$ in the augmentation ideal of $\h$,
$\tdelta\circ B^+(x)=B^+(x) \otimes \tun+(Id \otimes B^+)\circ \tdelta(x)$.
\end{enumerate}

{\bf Examples.}
\begin{eqnarray*}
\Delta(\tun)&=&\tun\otimes 1+1\otimes \tun,\\
\Delta(\tun\tun)&=&\tun\tun\otimes 1+1\otimes \tun\tun+\tun\otimes \tun,\\
\Delta(\tdeux)&=& \tdeux\otimes 1+1\otimes \tdeux+\tun \otimes \tun,\\
\Delta(\tdeux\tun)&=& \tdeux\tun\otimes 1+1\otimes \tdeux\tun+\tun \otimes \tun\tun+\tdeux\otimes \tun,\\
\Delta(\ttroisun)&=&\ttroisun\otimes 1+1\otimes \ttroisun+\tun\tun \otimes \tun+ \tun\otimes \tdeux,\\
\Delta(\ttroisdeux)&=&\ttroisdeux\otimes 1+1\otimes \ttroisdeux+\tdeux \otimes \tun+ \tun\otimes \tdeux.
\end{eqnarray*}

We now give a combinatorial description of this coproduct.
Let $F\in \F$ and $c \in Adm(F)$. Let $s_1 \geq_{h,l}\ldots\geq_{h,l} s_n$ be the vertices of $F$.
We shall say that $c$ is left-admissible if there exists $k \in \{1,\ldots,n\}$ such that the vertices of $P^c(F)$ are $s_1,\ldots,s_k$ and
the vertices of $R^c(F)$ are $s_{k+1},\ldots,s_n$. The set of left-admissible cuts of $F$ will be denoted  $Adm^l(F)$.

\begin{prop}
\label{10}
Let $F \in \F$. Then $\displaystyle \Delta(F)=\sum_{c \in {\cal A}dm^l(F)} P^c(F) \otimes R^c(F)+F \otimes 1+1 \otimes F$.
\end{prop}

{\bf Proof.} Consider $\Delta':\h \longrightarrow \h \otimes \h$, defined by the formula of the proposition \ref{10}. 
It is easy to show that, if $F,G \in \F$:
$$\left\{ \begin{array}{rcl}
\Delta(1)&=&1\otimes 1,\\
\Delta'(FG)&=&(F \otimes 1) \Delta'(G)+\Delta'(F)(1\otimes G)-F \otimes G,\\
\Delta'(B^+(F))&=&B^+(F) \otimes 1+(Id \otimes B^+)\circ \Delta'(F).
\end{array}\right.$$
By unicity in theorem \ref{9}, $\Delta=\Delta'$. $\Box$\\

Let us give a description of the branchs of the left-admissible cuts.

\begin{defi}
\textnormal{Let $F\in \F$. Let $I$ be a set of vertices of $F$. 
\begin{enumerate}
\item We shall say that $I$ is an ideal for $\geq_{high}$ if, for all vertices $s$, $s'$ of $F$:
$$ s\in I \mbox{ and } s'\geq_{high} s \: \Longrightarrow \: s'\in I.$$
\item We shall say that $I$ is an ideal for $\geq_{left}$ if, for all vertices $s$, $s'$ of $F$,
$$ s\in I \mbox{ and } s'\geq_{left} s \: \Longrightarrow \: s'\in I.$$
\item We shall say that $I$ is a biideal if $I$ is an ideal for $\geq_{high}$ and $\geq_{left}$.
\end{enumerate}}
\end{defi}

\begin{prop}
\label{12}
Let  $F \in \F$. Then $\displaystyle \Delta(F)=\sum_{\mbox{\scriptsize $I$ biideal of $F$}} I \otimes (F-I)$.
\end{prop}

{\bf Proof.} Similar to the proof of proposition \ref{10}. $\Box$\\

Let us precise the biideals of a forest $F$.

\begin{lemma}
Let $F\in \F$ and let $s_1 \geq_{h,l} \ldots \geq_{h,l} s_n$ be its vertices.
The biideals of $F$ are the sets $I_k=\{s_1,\ldots,s_k\}$, for $k\in \{0,\ldots, n\}$.
\end{lemma}

{\bf Proof.} Let $I$ be a biideal of $F$. Let $k$ be the greater integer such that $s_k \in I$. Then $I\subseteq I_k$.
Let $j \leq k$. Then $s_j \geq_{h,l} s_k$, so $s_j \geq_{high} s_k$ or $s_j \geq_{left} s_k$. As $I$ is a biideal, in both cases $s_j \in I$; hence, $I=I_k$.

it remains to show that $I_k$ is a biideal. Let $s_j \in I_k$ (so $j\leq k$), and $s_i$ be a vertex of $F$ such that $s_i \geq_{high} s_j$ or $s_i \geq_{left} s_j$.
Then, $s_i \geq_{h,l} s_j$ so $i\leq j \leq k$, and $s_i \in I_k$. $\Box$\\

{\bf Remark.} This implies that for any forest $F \in \F$, of weight $n$, for any $1<k<n$, 
there exists a unique left admissible cut $c$ such that the weight of $P^c(F)$ is equal to $k$.

\subsection{Antipode of $\h$}

As $\h$ is graded, with $\h_0=K$, it automatically has an antipode $S$, inductively defined by:
$$ \left\{ \begin{array}{rcl}
S(1)&=&1,\\
S(x)&=&-x-S(x')x'' \mbox{ if } \varepsilon(x)=0.
\end{array}
\right.$$

Because $\h$ is an infinitesimal Hopf algebra, $S$ satisfies $S(Ker(\varepsilon)^2)=0$, so $S(F)=0$  for all forest $F$ with at least two trees.
It remains to give a formula for the antipode of a single tree.

\begin{defi}
\textnormal{Let $t\in \T$. Let $s$ be the greatest vertex of $t$ for the total order relation $\geq_{h,l}$.
In other terms, $s$ is the leave of $t$ which is at most on the left.
\begin{enumerate}
\item Let $e$ be an edge of $t$. It will be called a {\it left edge} it it is on the path from the root to $s$.
\item Let $c$ be a (possibly empty) cut of $t$. We shall say that $c$ is a {\it left cut} if it cuts only left edges.
\end{enumerate}}
\end{defi}

Let $t \in \T$ and $c$ be a left cut of $t$. The cut $c$ makes $t$ into several trees $t_1,\ldots,t_n$.
These trees are indexed such that, by denoting $r_i$ the root of $t_i$ for all $i$, $r_1 \geq_{h,l} \ldots \geq_{h,l} r_n$ in $t$. 
The forest $t_1 \ldots t_n$ will be denoted $W^c(t)$. Moreover, we denote by $n_c$ the number of edges which are cut by $c$.

\begin{prop}
\label{15}
Let $t  \in \T$. Then $\displaystyle S(t)=- \sum_{\mbox{\scriptsize $c$ left cut of $t$}} (-1)^{n_c}\:  W^c(t)$.
\end{prop}

{\bf Proof.} We prove the result by induction on the weight $n$ of $t$. If $n=1$, then $t=\tun$ and the result is obvious.
Suppose the result true for all trees of weight $<n$. We put $t=B^+(t_1\ldots t_k)$. Two cases are possible.
\begin{enumerate}
\item If $k=1$, then:
\begin{eqnarray*}
\Delta(t)&=&t\otimes 1+(Id \otimes B^+)\circ \Delta(t_1)\\
&=&t \otimes 1+1 \otimes t+t_1 \otimes \tun+t_1' \otimes B^+(t_1''),\\
S(t)&=&-\underbrace{t}_{\substack{\mbox{\scriptsize $(-1)^{n_c} W^c(t)$,}\\ \mbox{\scriptsize $c$ empty }}}
-\underbrace{S(t_1)\tun}_{\substack{ \mbox{\scriptsize $(-1)^{n_c} W^c(t)$,}\\ \mbox{\scriptsize $c$ cuts the left edge}\\ \mbox{\scriptsize from the root}}}
-\underbrace{S(t'_1) B^+(t_1'').}_{\substack{\mbox{\scriptsize $(-1)^{n_c} W^c(t)$,}
\\ \mbox{\scriptsize $c$ does not cut the left}\\ \mbox{\scriptsize edge from the root}}}
\end{eqnarray*}
\item If $k\geq 2$, then:
\begin{eqnarray*}
\Delta(t)&=&t\otimes 1+(Id \otimes B^+)\circ \Delta(t_1\ldots t_k)\\
&=&t \otimes 1 +1 \otimes t+t_1 \ldots t_k \otimes \tun +(t_1\ldots t_k)' \otimes B^+((t_1\ldots t_k)'')\\
&=&t \otimes 1 +1 \otimes t+t_1 \ldots t_k \otimes \tun\\
&&+\sum_{i=1}^{k-1} t_1 \ldots t_i \otimes B^+(t_{i+1}\ldots t_k)+\sum_{i=1}^k t_1\ldots t_i' \otimes B^+(t_i'' \ldots t_k)\\
S(t)&=&-t-S(t_1\ldots t_k) \tun-\sum_{i=1}^{k-1} S(t_1 \ldots t_i) B^+(t_{i+1}\ldots t_k)-\sum_{i=1}^k S(t_1\ldots t_i') B^+(t_i'' \ldots t_k)\\
&=&-t-S(t_1)B^+(t_2 \ldots t_k)-S(t'_1)B^+(t_1''\ldots t_k)\\
&=&-\underbrace{t}_{\substack{\mbox{\scriptsize $(-1)^{n_c} W^c(t)$,}\\ \mbox{\scriptsize $c$ empty }}}
-\underbrace{S(t_1)B^+(t_2 \ldots t_k)}_{\substack{\mbox{\scriptsize $(-1)^{n_c} W^c(t)$,}\\ \mbox{\scriptsize $c$ cuts the left edge}
\\ \mbox{\scriptsize from the root}}}-\underbrace{S(t'_1) B^+(t_1''t_2 \ldots t_k).}_{\substack{\mbox{\scriptsize $(-1)^{n_c} W^c(t)$,}\\
\mbox{\scriptsize $c$ does not cut the left edge}\\ \mbox{\scriptsize from the root}}}
\end{eqnarray*}
\end{enumerate}
So the result holds for all forests. $\Box$\\

{\bf Examples.} 
\begin{eqnarray*}
S(\tun)&=&-\tun,\\
S(\tdeux)&=&-\tdeux+\tun\tun,\\
S(\ttroisun)&=&-\ttroisun+\tun\tdeux,\\
S(\ttroisdeux)&=&-\ttroisdeux+\tun\tdeux+\tdeux\tun-\tun\tun\tun,\\
S(\tquatreun)&=&-\tquatreun+\tun\ttroisun,\\
S(\tquatredeux)&=&-\tquatredeux+\tun\ttroisun+\tdeux\tdeux-\tun\tun\tdeux,\\
S(\tquatretrois)&=&-\tquatretrois+\tun\ttroisdeux,\\
S(\tquatrecinq)&=&-\tquatrecinq+\tun\ttroisdeux+\tdeux\tdeux+\ttroisdeux\tun
-\tun\tun\tdeux-\tun\tdeux\tun-\tun\tdeux\tun+\tun\tun\tun\tun.
\end{eqnarray*}

\subsection{Infinitesimal universal property}

\begin{theo}[Infinitesimal universal property]
\label{16}
Let $A$ be an infinitesimal Hopf algebra and let $L:A\longrightarrow A$ satisfying $\Delta(L(x))=L(x) \otimes 1+(Id \otimes L)\circ \Delta(x)$.
Then there exists a unique infinitesimal Hopf algebra morphism $\phi:\h \longrightarrow A$ such that $\phi \circ B^+=L\circ \phi$. 
\end{theo}

{\bf Proof.} By the universal property of $\h$, there exists a unique algebra morphism $\phi$ satisfying $\phi \circ B^+=L\circ \phi$. 
Let us show that it is also a coalgebra morphism.\\

{\it First step.} We first show that $\varepsilon \circ L=0$. Let $a \in A$.
\begin{eqnarray*}
\varepsilon \circ L(a)&=&(\varepsilon \otimes \varepsilon)\circ \Delta \circ L(a)\\
&=&(\varepsilon \otimes \varepsilon)(L(a) \otimes 1+a^{(1)} \otimes L(a^{(2)}))\\
&=&\varepsilon\circ L(a)+\varepsilon(a^{(1)}) \varepsilon\circ L(a^{(2)})\\
&=&\varepsilon\circ L(a)+\varepsilon\circ L(a).
\end{eqnarray*}
So $\varepsilon \circ L(a)=0$.\\

{\it Second step.} We consider $X=\{x\in \h\:/\: \varepsilon \circ \phi(a)=\varepsilon(x)\}$.
As $\varepsilon \circ \phi$ and $\varepsilon$ are both algebra morphisms, $X$ is a subalgebra of $\h$. Let $x \in \h$. Then:
$$(\varepsilon \circ \phi)(B^+(x))=\varepsilon \circ L \circ \phi(x)=0=\varepsilon(B^+(x)).$$
So $Im(B^+)\subseteq X$. As $X$ is a subalgebra, $X=\h$, and $\varepsilon \circ \phi=\varepsilon$.\\

{\it Third step.} We consider $Y=\{x \in \h\:/\: \Delta \circ \phi(x)=(\phi \otimes \phi)\circ \Delta\}$.
As $\Delta \circ \phi$ and $(\phi \otimes \phi)\circ \Delta$ are algebra morphisms from $\h$  to $(A\otimes A,.\varepsilon)$, 
$Y$ is a subalgebra of $\h$. Let $x \in Y$. 
\begin{eqnarray*}
(\Delta \circ \phi)(B^+(x))&=&\Delta \circ L \circ \phi(x)\\
&=&L\circ \phi(x) \otimes 1+\phi(x)^{(1)} \otimes L(\phi(x)^{(2)})\\
&=&\phi \circ B^+(x) \otimes 1+\phi(x^{(1)}) \otimes L(\phi(x^{(2)}))\\
&=&\phi \circ B^+(x) \otimes \phi(1)+\phi(x^{(1)}) \otimes \phi(B^+(x^{(2)}))\\
&=&(\phi \otimes \phi)(B^+(x) \otimes 1+x^{(1)}\otimes B^+(x^{(2)}))\\
&=&(\phi \otimes \phi) \circ \Delta (B^+(x)).
\end{eqnarray*}
So $B^+(x) \in Y$. As $Y$ is a subalgebra of $\h$ stable under $B^+$, $Y=\h$. Hence, $\phi$ is a coalgebra morphism. $\Box$

\subsection{A pairing on $\h$}

\begin{defi}
\textnormal{ The application $\gamma$ is defined by:
$$\gamma: \left\{ \begin{array}{rcl}
\h&\longrightarrow & \h\\
t_1\ldots t_n \in \F&\longrightarrow & \delta_{t_1,\tun} t_2\ldots t_n.
\end{array}\right.$$}
\end{defi}

\begin{lemma}
\label{18}
\begin{enumerate}
\item $\gamma$ is homogeneous of degree $-1$.
\item For all $x,y\in \h$, $\gamma(xy)=\gamma(x)y+\varepsilon(x) \gamma(y)$.
\item $Ker(\gamma) \cap Prim(\h)=(0)$.
\end{enumerate}
\end{lemma}

{\bf Proof.} \begin{enumerate}
\item Trivial.
\item Immediate for $x,y\in \F$, separating the cases $x=1$ and $x\neq 1$.
\item Let us take $p\in Ker(\gamma)$, non-zero, and primitive. Then $p$ can be written as:
$$p=\sum_{F\in \F} a_F F.$$
Let us choose a forest $F=t_1\ldots t_n$ such that:
\begin{enumerate}
\item $a_F \neq 0$.
\item If $G=t'_1\ldots t'_m \in \F$ is such that $a_G \neq 0$, then $m\leq n$. If moreover $m=n$,
we put $t_1=B^+(s_1\ldots s_k)$ and $t'_1=B^+(s'_1\ldots s'_l)$; then $k \leq l$.
\end{enumerate}
Suppose that $t_1\neq \tun$. Then $k\neq 0$. We consider the cut $c$ on the edge from the root of $t_1$ to the root of $s_1$.
This is a left-admissible cut, so $s_1 \otimes B^+(s_2 \ldots s_m)t_2\ldots t_n$ appears in $\Delta(F)$.
Because $p$ is primitive, there exists another forest $G$ such that $a_G\neq 0$ 
and $s_1 \otimes B^+(s_2\ldots s_m) t_2 \ldots t_n$ appears in $\Delta(G)$. Three cases are possible:
\begin{enumerate}
\item $G=s_1 B^+(s_2\ldots s_m) t_2 \ldots t_n$: contradicts the maximality of $n$ for $F$.
\item $G=F$: contradicts that $G \neq F$.
\item $G$ is obtained by grafting $s_1$ on a vertex of $s_2$: contradicts the minimality of $k$.
\end{enumerate}
In every case, we obtain a contradiction. So $t_1=\tun$. Hence, $\gamma(p)\neq 0$. $\Box$
\end{enumerate}

\begin{theo}
\label{19}
There exists a unique pairing $\langle-,-\rangle:\h \times \h \longrightarrow K$, satisfying:
\begin{description}
\item[\it 1.] $\langle1,x\rangle =\varepsilon(x)$ for all $x \in \h$.
\item[\it 2.] $\langle xy,z\rangle =\langle y\otimes x, \Delta(z)\rangle $ for all $x,y,z \in \h$.
\item[\it 3.] $\langle B^+(x),y\rangle =\langle x,\gamma(y)\rangle $ for all $x,y \in \h$.
\end{description}
Moreover:
\begin{description}
\item[\it 4.] $\langle-,-\rangle $ is symmetric and non-degenerate.
\item[\it 5.] If $x$ and $y$ are homogeneous of different weights, $\langle x,y\rangle =0$.
\item[\it 6.] $\langle S(x),y\rangle =\langle x,S(y)\rangle $ for all $x,y \in \h$.
\end{description}
\end{theo}

{\bf Proof.} 

{\it Unicity.} Assertions 1-3 entirely determine $\langle F,G\rangle $ for $F,G \in \F$ by induction on the weight.\\

{\it Existence.} We consider the graded infinitesimal Hopf algebra $A=\h^{*,op,cop}$.
As $\gamma$ is homogeneous of degree $-1$, it can be transposed in an application:
$L=\gamma^*:A\longrightarrow A$. This linear application is homogeneous of degree $1$. Let $f\in A^*$, and $x,y \in \h$.
\begin{eqnarray*}
\left((\Delta \circ L)(f)\right)(x \otimes y)&=&L(f)(yx)\\
&=&f(\gamma(yx))\\
&=&f(\gamma(y)x +\varepsilon(y)\gamma(x))\\
&=&\Delta(f)(x \otimes \gamma(y))+ (f \otimes 1)(\gamma(x) \otimes y)\\
&=&\left((Id \otimes L)\circ \Delta(f)+L(f) \otimes 1\right)(x \otimes y).
\end{eqnarray*}
So $\Delta\circ L(f)=(Id \otimes L)\circ \Delta(f)+L(f) \otimes 1$ for all $f \in A$.
By the infinitesimal universal property, there exists a unique infinitesimal Hopf algebra morphism $\phi:\h \longrightarrow A$, 
such that $\phi \circ B^+=L\circ \phi$. We then put $\langle x,y\rangle =\phi(x)(y)$ for all $x,y \in \h$. Let us show that this pairing satisfies 1-6.
\begin{description}
\item[\textnormal{1.}] For all $x \in \h$, $\langle 1,x\rangle =\phi(1)(x)=\varepsilon(x)$.
\item[\textnormal{2.}]  For all $x,y,z \in \h$:
$$\langle xy,z\rangle =\phi(xy)(z)=(\phi(x)\phi(y))(z)=(\phi(y) \otimes \phi(x))( \Delta(z))=\langle y\otimes x,\Delta(z)\rangle.$$
\item[\textnormal{3.}] For all $x,y \in \h$:
$$\langle B^+(x),y\rangle =(\phi \circ B^+(x))(y)=(L\circ \phi(x))(y)=\phi(x)(\gamma(y))=\langle x,\gamma(y)\rangle.$$
\item[\textnormal{5.}] As $L$ is homogeneous of degree $1$, $\phi$ is homogeneous of degree $0$. This implies 5.
\item[\textnormal{6.}] As $\phi$ is an infinitesimal Hopf algebra morphism and $S^*$ is the antipode of $A$, $S^* \circ \phi=\phi \circ S$. 
Hence, for all $x,y\in \h$:
$$\langle S(x),y\rangle =\phi(S(x))(y)=(S^* \circ \phi(x))(y) =\phi(x)(S(y))=\langle x,S(y)\rangle.$$
\item[\textnormal{4.}] Let us first show that $\langle-,-\rangle $ is symmetric.
We put $\langle x,y\rangle '=\langle y,x\rangle $ for all $x,y \in \h$. Let us show that $\langle-,-\rangle '$ satisfies 1-3. 
For all $x \in \h$, as $f(1)=\varepsilon(f)$ for all $f \in \h^*$:
$$\langle 1,x\rangle '=\langle x,1\rangle =\phi(x)(1)=\varepsilon(\phi(x))=\varepsilon(x).$$
For all $x,y,z \in \h$, as $\phi$ is a coalgebra morphism:
\begin{eqnarray*}
\langle xy,z\rangle '&=&\langle z,xy\rangle \\
&=&\phi(z)(xy)\\
&=&(\Delta\circ \phi(z))(y \otimes x)\\
&=&((\phi \otimes \phi)\circ \Delta(z))(y \otimes x)\\
&=&\langle \Delta(z),y\otimes x\rangle \\
&=&\langle y \otimes x,\Delta(z)\rangle '.
\end{eqnarray*}

Let us show 3 for $\langle-,-\rangle '$ with $y \in \F$, by induction on $n=weight(y)$.
If $n=0$, then $y=1$:
$$\langle B^+(x),y\rangle '=\langle 1,B^+(x)\rangle=\varepsilon\circ B^+(x)=0=\langle x,\gamma(y)\rangle'.$$
Suppose the result true for every forest $F$ of weight $<n$. Two cases are possible:
\begin{description}
\item[\textnormal{-}]  $y=B^+(z)$. We can restrict to the case where $x$ is also a forest. Then:
$$\langle B^+(x),y\rangle '=\langle B^+(z),B^+(x)\rangle =\langle z,\gamma\circ B^+(x)\rangle=\delta_{x,1}\langle z,1\rangle =\varepsilon(x) \varepsilon(z).$$
Moreover, $\langle x,\gamma(y)\rangle '=\langle \gamma \circ B^+(z),x\rangle=\delta_{z,1}\langle 1,x\rangle=\varepsilon(z) \varepsilon(x)$.
\item[\textnormal{-}]  $y$ is a forest with at least two trees. 
Then $y$ can be written $y=y_1y_2$, with the induction hypothesis avalaible for $y_1$ and $y_2$. Then:
\begin{eqnarray*}
\langle B^+(x),y\rangle '&=&\langle y_1y_2,B^+(x)\rangle \\
&=&\langle y_2 \otimes y_1,\Delta \circ B^+(x)\rangle \\
&=&\langle y_2 \otimes y_1,B^+(x) \otimes 1+(Id \otimes B^+)\circ \Delta(x)\rangle \\
&=&\langle \gamma(y_2) \varepsilon(y_1),x\rangle +\langle y_2 \otimes \gamma(y_1),\Delta(x)\rangle \\
&=&\langle \varepsilon(y_1)\gamma(y_2)+\gamma(y_1)y_2,x\rangle \\
&=&\langle \gamma(y),x\rangle \\
&=&\langle x,\gamma(y)\rangle '.
\end{eqnarray*}
\end{description}
So $\langle-,-\rangle'$ satisfies 1-3. By unicity, $\langle-,-\rangle '=\langle-,-\rangle$, so $\langle-,-\rangle$ is symmetric.

Note that this implies that for all $x,y \in \h$:
$$(\phi\circ \gamma(x))(y)=\langle \gamma(x),y\rangle =\langle y,\gamma(x)\rangle =\langle B^+(y),x\rangle=\langle x,B^+(y)\rangle =(\phi(x))(B^+(y)).$$
So $\phi\circ \gamma=(B^+)^* \circ \phi$.

It remains to prove that $\langle-,-\rangle $ is non degenerate. It is equivalent to show that $\phi$ is monic. 
Suppose that it is not. Let us choose a non-zero element $p\in Ker(\phi)$ of lowest degree. 
As $\phi$ is a coalgebra morphism, its kernel is a coideal, so $p \in Prim(\h)$.
By lemma \ref{18}-3, $\gamma(p) \neq 0$. Moreover, $\phi\circ \gamma(p)=(B^+)^* \circ \phi(p)=0$. 
So $\gamma(p) \in Ker(\phi)$, is non-zero, of degree strictly smaller than $p$: 
this contradicts the choice of $p$. So $\phi$ is monic. $\Box$
\end{description}

{\bf Remark.} Similarly with the usual case, it is possible to define a pairing between $\h$ and itself, using the application:
$$\gamma': \left\{ \begin{array}{rcl}
\h&\longrightarrow & \h\\
t_1\ldots t_n \in \F&\longrightarrow & \delta_{t_n,\tun} t_1\ldots t_{n-1}.
\end{array}\right.$$
Unhappily, this pairing is degenerate: it is for example not difficult to show that the primitive element $\ttroisun-\tun\tdeux$ belongs to $\h^\perp$.

\begin{defi}
\textnormal{We denote by $(f_F)_{F\in \F}$ the dual basis of the basis of forests. In other terms, for all $F \in \F$, 
$f_F$ is defined by $\langle f_F,G\rangle =\delta_{F,G}$, for all forest $G \in \F$.}
\end{defi}

\begin{prop}
\label{21}
\begin{enumerate}
\item For all forest $F\in \F$, $B^+(f_F)=f_{\tun F}$.
\item For all forest $F\in \F$:
$$\gamma(f_F)= \left\{ \begin{array}{rcl}
0&\mbox{if}&F\notin \T,\\
f_{B^-(F)}&\mbox{if}&F\in \T,
\end{array}\right.$$
where $B^-(F)$ is the forest obtained by deleting the root of $F$.
\item For all forest $F \in \F$, $\displaystyle \Delta(f_F)=\sum_{\substack{F_1,F_2 \in \F\\F_1F_2=F}}f_{F_2} \otimes f_{F_1}$.
\end{enumerate}
\end{prop}

{\bf Proof.} \begin{enumerate}
\item Let $G \in \F$. Then:
\begin{eqnarray*}
\langle B^+(f_F),G\rangle &=&\langle f_F,\gamma(G)\rangle \\
&=&\left\{\begin{array}{l}
0 \mbox{ if $G$ is not of the form $\tun H$},\\
\delta_{F,H} \mbox{ if $G=\tun H$},
\end{array}\right.\\
&=&\delta_{\tun F,G}\\
&=&\langle f_{\tun F},G\rangle.
\end{eqnarray*}
As $\langle-,-\rangle$ is non-degenerate, $B^+(f_F)=f_{\tun F}$. 
\item Suppose first that $F$ is not a tree. Then, for all $G \in \F$:
$$\langle \gamma(f_F),G\rangle =\langle f_F,B^+(G)\rangle =\delta_{F,B^+(G)}=0.$$
So $\gamma(f_F)=0$. Suppose now that $F$ is a tree. Then, for all $G \in \F$:
$$\langle \gamma(f_F),G\rangle =\delta_{F,B^+(G)}=\delta_{B^-(F),G}=\langle f_{B^-(F)},G\rangle.$$
So $\gamma(f_F)=f_{B^-(F)}$. 
\item Indeed, for all forests $G_1,G_2\in \F$,
$$\langle \Delta(f_F),G_1 \otimes G_2\rangle =\langle f_F,G_2G_1\rangle =\delta_{F,G_2G_1}
=\sum_{\substack{F_1,F_2 \in \F\\F_1F_2=F}}\langle f_{F_2} \otimes f_{F_1},G_1 \otimes G_2\rangle.$$
As $\langle-,-\rangle $ is non-degenerate, this proves the last point. $\Box$
\end{enumerate}

\begin{prop}
The familly $(f_t)_{t\in \T}$ is a basis of $Prim(\h)$.
\end{prop}

{\bf Proof.} Immediate corollary of proposition \ref{21}-3. $\Box$\\

As an example of decomposition in the dual basis, we give the following result:

\begin{cor}
\label{23} 
For all $n \in \mathbb{N}$, $\displaystyle \tun^n=\sum_{\substack{F\in \F\\ weight(F)=n}}f_F$.
\end{cor}

{\bf Proof.} For all $n \in \mathbb{N}^*$,  we can put, by homogeneity:
$$\tun^n=\sum_{\substack{F\in \F\\ weight(F)=n}}a_F f_F.$$
We define $X=\{F\in \F\:/\: a_F=1\}$. Note that $f_1=1$, so $1\in X$. Let $F_1,F_2 \in X$. Then, if $n=weight(F_1F_2)$:
\begin{eqnarray*}
a_{F_1F_2}&=&\langle \tun^n, F_1F_2 \rangle\\
&=&\langle \Delta(\tun^n), F_2 \otimes F_1 \rangle\\
&=&\langle \sum_{i=0}^n \tun^i \otimes \tun^{n-i}, F_2 \otimes F_1 \rangle\\
&=&\langle \tun^{weight(F_2)} \otimes \tun^{weight(F_1)},F_2\otimes F_1 \rangle +0\\
&=&a_{F_2}a_{F_1}\\
&=&1.
\end{eqnarray*}
So $F_1F_2 \in X$. Moreover, if $m=weight(F_1)$:
$$a_{B^+(F_1)}=\langle \tun^{m+1}, B^+(F_1) \rangle=\langle \gamma(\tun^{m+1}),F_1 \rangle=\langle \tun^m,F_1 \rangle=a_{F_1}=1.$$
So $B^+(F_1)\in X$. Hence, $1 \in X$, and $X$ is stable by product and by $B^+$. So $X=\F$. $\Box$

\subsection{Combinatorial interpretation of the pairing}

{\bf Notation.} Let $F \in \F$. We denote by $Vert(F)$ the set of vertices of $F$.

\begin{theo}
\label{24}
Let $F,G\in \F$. Let $S(F,G)$ be the set of bijections $\sigma:Vert(F) \longrightarrow Vert(G)$ such that, for all vertices $x$, $y$ of $F$:
\begin{enumerate}
\item $(x \leq_{high} y)$ $\Longrightarrow$ $(\sigma(x) \geq_{left} \sigma(y))$.
\item $(x \leq_{left} y)$ $\Longrightarrow$ $(\sigma(x) \geq_{h,l} \sigma(y))$.
\item $(\sigma(x) \leq_{high} \sigma(y))$ $\Longrightarrow$ $(x \geq_{left} y)$.
\item $(\sigma(x) \leq_{left} \sigma(y))$ $\Longrightarrow$ $(x \geq_{h,l} y)$.
\end{enumerate}
Then $\langle F,G \rangle=card(S(F,G))$.
\end{theo}

{\bf Proof.} If $F$ and $G$ have different weights, as $\langle-,-\rangle$ is homogeneous, then $\langle F,G \rangle=0$ and $S(F,G)$ is empty, 
so the result holds. Let us suppose now that $F$ and $G$ have the same weight $n$ and let us proceed by induction on $n$.
If $n=0$, then $F=G=1$, and the result holds. For the hereditary, we have the following cases.
\begin{enumerate}
\item $F=B^+(F_1)$. We have the two following subcases.
\begin{enumerate}
\item $G$ is not of the form $\tun G_1$. Then $\gamma(G)=0$, so 
$\langle F,G \rangle=\langle B^+(F_1),G \rangle=\langle F_1,\gamma(G) \rangle=0$.
Let us assume that $S(F,G)$ is not empty, and let $\sigma \in S(F,G)$. Let $r$ be the root of $F$.
For all $x \in vert(F)$, $x \geq_{high} r$, so $\sigma(r) \geq_{left} \sigma(x)$. Hence, as $\sigma$ is epic, $G$ is of the form $\tun G_1$: contradiction.
So $S(F,G)=\emptyset$, and the result holds.
\item $G=\tun G_1$. Then $\langle F,G\rangle=\langle F_1,G_1 \rangle$.
Let $\sigma \in S(F,G)$. As in the preceding point, $\sigma(r)$ is the vertex of $\tun$, so we can consider the application:
$$\Psi: \left\{ \begin{array}{rcl}
S(F,G)&\longrightarrow & S(F_1,G_1)\\
\sigma&\longrightarrow &\sigma_{\mid Vert(F_1)}.
\end{array}\right.$$
It is obviously monic. Let us show it is epic. Let $\sigma_1 \in S(F_1,G_1)$, and let $\sigma:Vert(F)\longrightarrow Vert(G)$
extending $\sigma_1$ by sending the root of $F$ to the vertex of $\tun$. Let us show that $\sigma \in S(F,G)$. Let $x,y \in Vert(F)$.
We can suppose they are distinct.
\begin{description}
\item[\textnormal{-}] If $x \leq_{high} y$, two cases are possible. If $x$ is the root of $F$, then $\sigma(x) \geq_{left} \sigma(y)$.
If not, then $x$ and $y$ are vertices of $F_1$, so $\sigma(x) \geq_{left} \sigma(y)$ in $G_1$, hence in $G$.
\item[\textnormal{-}] If $x \leq_{left} y$, then both of them are vertices of $F_1$, so $\sigma(x) \geq_{h,l} \sigma(y)$ in $G_1$, hence in $G$.
\item[\textnormal{-}] If $\sigma(x) \leq_{high} \sigma(y)$, then $\sigma(x)$ and $\sigma(y)$ are vertices of $G_1$,
so $x \geq_{left} y$ in $F_1$, hence in $F$.
\item[\textnormal{-}] If $\sigma(x) \leq_{left} \sigma(y)$, then two cases are possible.
If $\sigma(y)$ is the vertex of $\tun$, then $y$ is the root of $F$, so $x \geq_{h,l} y$. If not, then $\sigma(x)$ and $\sigma(y)$ are vertices of $G_1$,
so $x \geq_{h,l} y$ in $F_1$, hence in $F$.
\end{description}
So $\sigma \in S(F,G)$, and $\Psi(\sigma)=\sigma_1$. So $\Psi$ is a bijection. As a consequence:
$$card(S(F,G))=card(S(F_1,G_1))=\langle F_1,G_1\rangle=\langle F,G \rangle.$$
\end{enumerate}
\item $F=F_1F_2$, with $F_1,F_2 \neq 1$. Let $J$ be the unique biideal of $G$ with the same weight as $F_2$. Then:
$$\langle F,G \rangle=\langle F_2 \otimes F_1, \sum_{\mbox{\scriptsize{$I$ biideal of $G$}}} I \otimes (G-I)\rangle
=\langle F_2,J \rangle \langle F_1,G-J\rangle.$$
Let $\sigma \in S(F,G)$. Let $\sigma(x) \in \sigma(Vert(F_2))$, and $y' \in vert(G)$, such that $y' \geq_{high} \sigma(x)$ or $y' \geq_{left} \sigma(x)$.
As $\sigma$ is epic, we put $y'=\sigma(y)$. Then $x \geq_{high} y$ or $x \geq_{h,l} y$. In both cases, as $x \in Vert(F_2)$, $y\in Vert(F_2)$.
So $ \sigma(Vert(F_2))$ is a biideal of $G$. Considering its weight, it is $J$. We can consider the application:
$$\Phi: \left\{ \begin{array}{rcl}
S(F,G)&\longrightarrow & S(F_1,G-J) \times S(F_2,J)\\
\sigma &\longrightarrow &(\sigma_{\mid vert(F_1)},\sigma_{\mid vert(F_2)}).
\end{array}\right.$$
It is clearly monic. Let be $(\sigma_1,\sigma_2) \in S(F_1,G-J) \times S(F_2,J)$ and let $\sigma:Vert(F) \longrightarrow Vert(G)$,
such that $\sigma_{\mid vert(F_i)}=\sigma_i$ for $i=1,2$. Let us show that $\sigma\in S(F,G)$. Let $x,y \in Vert(F)$.
\begin{description}
\item[\textnormal{-}] If $x \leq_{high} y$, then $x,y \in F_1$ or $x,y \in F_2$. So $x\geq_{left} y$ in $J$ or in $G-J$, hence in $G$.
\item[\textnormal{-}] If $x \leq_{left} y$, two cases are possible. If $x,y \in F_1$ or $x,y \in F_2$, then $x\geq_{h,l} y$ in $J$ or in $G-J$, hence in $G$.
If $x \in Vert(F_2)$ and $y \in Vert(F_1)$, then $\sigma(x) \in J$ and $\sigma(y) \in G-J$. 
As $J$ is a biideal, $\sigma(y) \geq_{h,l} \sigma(x)$ is impossible. As $\geq_{h,l}$ is a total order, $\sigma(x) \geq_{h,l} \sigma(y)$.
\item[\textnormal{-}] If $\sigma(x) \leq_{high} \sigma(y)$, two cases are possible. If $\sigma(x), \sigma(y) \in J$ or
 $\sigma(x), \sigma(y) \in G-J$, then $x \geq_{left} y$ in $F_1$ or in $F_2$, hence in $F$. 
If $\sigma(x) \in G-J$ and $\sigma(y) \in J$, then $x \in vert(F_1)$ and $y\in vert(F_2)$, so $x \geq_{left} y$.
\item[\textnormal{-}] If $\sigma(x) \leq_{left} \sigma(y)$,  two cases are possible. If $\sigma(x), \sigma(y) \in J$ or
 $\sigma(x), \sigma(y) \in G-J$, then $x \geq_{h,l} y$ in $F_1$ or in $F_2$, hence in $F$. 
 If $\sigma(x) \in G-J$ and $\sigma(y) \in J$, then $x \in vert(F_1)$ and $y\in vert(F_2)$, so $x \geq_{h,l} y$.
\end{description}
So $\sigma \in S(F,G)$, and $\Phi(\sigma)=(\sigma_1,\sigma_2)$. So $\Phi$ is a bijection. Hence:
$$\langle F,G \rangle=\langle F_2,J \rangle \langle F_1,G-J\rangle=card(S(F_1,G-J))card(S(F_2,J))=card(S(F,G)).$$
\end{enumerate}
So the result holds for all $n$. $\Box$\\

{\bf Remarks.}
\begin{enumerate}
\item There is obviously a bijection:
$$ \left\{ \begin{array}{rcl}
S(F,G)&\longrightarrow & S(G,F)\\
\sigma&\longrightarrow &\sigma^{-1}.
\end{array}\right.$$
This gives another proof of the symmetry of $\langle-,-\rangle$.
\item Let us assume that $S(F,G)$ is not empty, and let $\sigma\in S(F,G)$. By definition, if $x \leq_{h,l} y$ in $Vert(F)$,
then $\sigma(x) \geq_{h,l} \sigma(y)$ in $Vert(G)$, so $\sigma$ is the unique decreasing bijection from $(Vert(F),\geq_{h,l})$
to  $(Vert(G),\geq_{h,l})$. So, for any forests $F,G \in \F$, $\langle F,G\rangle=0$ or $1$.
\end{enumerate}

\section{Poset of forests and applications}

\subsection{Partial order on $\F$}

{\bf Notations.} Let $n \in \mathbb{N}$. The set of planar forests with $n$ vertices will be denoted by $\F(n)$.

\begin{defi} 
\label{25} \textnormal{Let $F \in \F$. 
\begin{enumerate}
\item An {\it admissible transformation} on $F$ is a local transformation of one of the following type 
(the part of $F$ which is not in the frame remains unchanged):
$$\begin{array}{rccc}
\mbox{First kind: }&
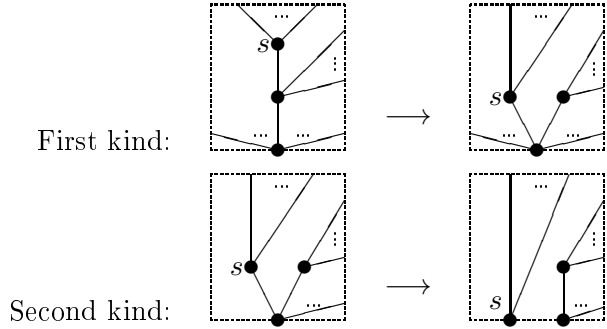
\begin{picture}(55,60)(-25,0)
\put(0,0){\circle*{5}}
\put(0,0){\line(0,1){20}}
\put(0,20){\circle*{5}}
\put(0,20){\line(0,1){20}}
\put(0,40){\circle*{5}}
\put(0,40){\line(-1,1){15}}
\put(0,40){\line(1,1){15}}
\put(-8,37){$s$}
\put(-2,50){.}
\put(0,50){.}
\put(2,50){.}
\put(0,0){\line(-4,1){25}}
\put(0,0){\line(4,1){25}}
\put(-10,5){.}
\put(-8,5){.}
\put(-6,5){.}
\put(10,5){.}
\put(8,5){.}
\put(6,5){.}
\put(0,20){\line(4,1){25}}
\put(0,20){\line(1,1){25}}
\put(20,30){.}
\put(20,32){.}
\put(20,34){.}
\put(-25,0){\dashbox{1}(50,55)}
\end{picture}
&\begin{picture}(22,0)(0,0)
\put(0,10){$ \longrightarrow $}
\end{picture}
&\begin{picture}(55,60)(-25,0)
\put(0,0){\circle*{5}}
\put(0,0){\line(1,2){10}}
\put(0,0){\line(-1,2){10}}
\put(-10,20){\circle*{5}}
\put(10,20){\circle*{5}}
\put(-18,17){$s$}
\put(0,0){\line(-4,1){25}}
\put(0,0){\line(4,1){25}}
\put(-14,5){.}
\put(-12,5){.}
\put(-10,5){.}
\put(14,5){.}
\put(12,5){.}
\put(10,5){.}
\put(-10,20){\line(0,1){35}}
\put(-10,20){\line(2,3){23.5}}
\put(-2,50){.}
\put(0,50){.}
\put(2,50){.}
\put(10,20){\line(3,5){15}}
\put(10,20){\line(4,1){15}}
\put(20,28){.}
\put(20,30){.}
\put(20,32){.}
\put(-25,0){\dashbox{1}(50,55)}
\end{picture}\\
\mbox{Second kind: }&
\begin{picture}(55,60)(-25,0)
\put(0,0){\circle*{5}}
\put(0,0){\line(1,2){10}}
\put(0,0){\line(-1,2){10}}
\put(-10,20){\circle*{5}}
\put(10,20){\circle*{5}}
\put(-18,17){$s$}
\put(0,0){\line(4,1){25}}
\put(14,5){.}
\put(12,5){.}
\put(10,5){.}
\put(-10,20){\line(0,1){35}}
\put(-10,20){\line(2,3){23.5}}
\put(-2,50){.}
\put(0,50){.}
\put(2,50){.}
\put(10,20){\line(3,5){15}}
\put(10,20){\line(4,1){15}}
\put(20,28){.}
\put(20,30){.}
\put(20,32){.}
\put(-25,0){\dashbox{1}(50,55)}
\end{picture}
&\begin{picture}(22,0)(0,0)
\put(0,10){$ \longrightarrow $}
\end{picture}
&\begin{picture}(55,60)(-25,0)
\put(10,0){\circle*{5}}
\put(10,20){\circle*{5}}
\put(-10,0){\circle*{5}}
\put(-18,4){$s$}
\put(10,0){\line(0,1){20}}
\put(10,20){\line(3,5){15}}
\put(10,20){\line(4,1){15}}
\put(20,28){.}
\put(20,30){.}
\put(20,32){.}
\put(-10,0){\line(0,1){55}}
\put(-10,0){\line(2,5){22}}
\put(-2,50){.}
\put(0,50){.}
\put(2,50){.}
\put(10,0){\line(3,1){15}}
\put(-25,0){\dashbox{1}(50,55)}
\put(13,5){.}
\put(15,5){.}
\put(17,5){.}
\end{picture}
\end{array}$$
Such a transformation will be said to hold on the vertex $s$.
\item Let $s_1 \geq_{h,l} \ldots \geq_{h,l} s_n$ be the vertices of $F$. An admissible transformation on $F$
will be said to be an $i$-transformation if it holds on the vertex $s_i$.
\end{enumerate}}
\end{defi}

{\bf Example.} Let $t=\tquatredeux$. The $1$-transformation turns $t$ into $\tquatreun$.
The $2$-transformation turns $t$ into $\tdeux \tdeux$. There are no $3$- and $4$-transformation.

\begin{defi}
\textnormal{Let $I \subseteq \mathbb{N}^*$. Let $F,G \in \F$.
We shall say that $F\leq_I G$ if there exists a finite sequence $F_0,\ldots,F_k$ of elements of $\F$ such that:
\begin{enumerate}
\item For all $i \in \{0,\ldots, k-1\}$, $F_{i+1}$ is obtained from $F_i$ by a $j$-transformation, for a certain $j\in I$.
\item $F_0=F$.
\item $F_k=G$.
\end{enumerate}}
\end{defi}

\begin{prop}
For all $I\subseteq \mathbb{N}^*$, $\leq_I$ is a partial order $\F$.
\end{prop}

{\bf Proof.} Indeed, $\leq_I$ is transitive and reflexive. Let $F,G$ be two forests, such that $F \leq_I G$ and $G \leq_I F$. We put:
$$ \left\{ \begin{array}{rcl}
n_F&=&\displaystyle \sum_{\mbox{\scriptsize $s$ vertex of $F$}} height(s),\\
n_G&=&\displaystyle \sum_{\mbox{\scriptsize $s$ vertex of $G$}} height(s).
\end{array}\right.$$
As $F \leq_I G$, there exists $F_0,\ldots,F_k \in \F$, such that $F_0=F$, $F_k=G$, and $F_{i+1}$ is obtained from $F_i$ by an admissible transformation.
Each admissible transformation decreases the height of a vertex by $1$, so $n_G=n_F-k$, so $n_G\leq n_F$. 
As $G \leq_I F$, in the same way, $n_F \leq n_G$, so $n_F=n_G$, and $k=0$. As a consequence, $F=G$. $\Box$\\

{\bf Remarks.} \begin{enumerate}
\item $F$ and $G$ are comparable for $\leq_\emptyset$ if, and only if, they are equal.
\item We shall denote $\leq$ instead of $\leq_{\mathbb{N}}$. This order $\leq$ is the order generated by all the admissible transformations.
\item If $F$ and $G \in \F$ are comparable for one of these orders, they have the same weight.
So $\F= \displaystyle \bigcup_{n\in \mathbb{N}} \F(n)$ as a poset.\\
\end{enumerate}

{\bf Examples.}  The posets $\F(0)$ and $\F(1)$ are reduced to a single element. 
Here are the Hasse graphs of $(\F(2),\leq)$, $(\F(3),\leq)$ and $(\F(4),\leq)$:
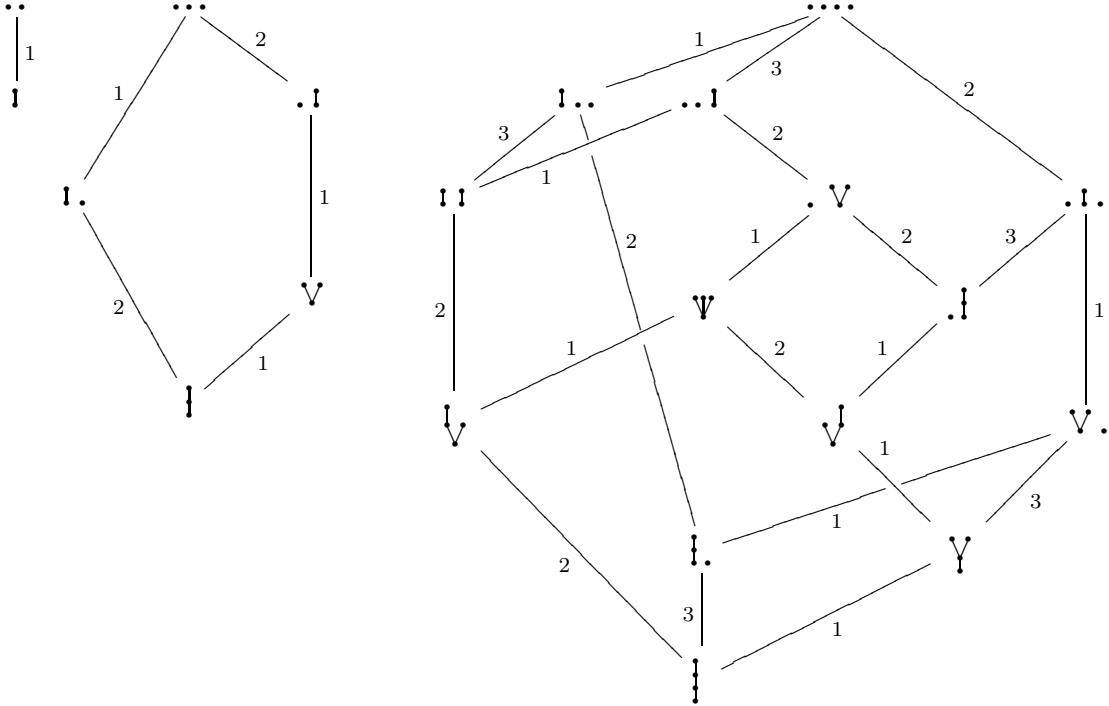
$$\xymatrix{\tun\tun \ar@{-}[d]^1\\\tdeux}\:
\xymatrix{&\tun\tun\tun \ar@{-}[ldd]_1 \ar@{-}[rd]^2&\\
&&\tun\tdeux\ar@{-}[dd]^1 \\
\tdeux\tun\ar@{-}[rdd]_2 &&\\
&&\ttroisun\ar@{-}[ld]^1 \\
&\ttroisdeux&}\:
\xymatrix{&&&&\tun\tun\tun\tun\ar@{-}[lld]_1\ar@{-}[ld]^3\ar@{-}[rrdd]^2&&\\
&&\tdeux\tun\tun\ar@{-}[ld]_3\ar@{-}|(.11)\hole|(.53)\hole[rdddd]^<(.3){2}
&\tun\tun\tdeux\ar@{-}[lld]^<(.6){1}\ar@{-}[rd]^2&&&\\
&\tdeux\tdeux\ar@{-}[dd]_2&&&\tun\ttroisun\ar@{-}[ld]_1 \ar@{-}[rd]^2&&
\tun\tdeux\tun\ar@{-}[ld]_3\ar@{-}[dd]^1\\
&&&\tquatreun\ar@{-}[lld]_1\ar@{-}[rd]^2&&\tun\ttroisdeux\ar@{-}[ld]_1&\\
&\tquatredeux\ar@{-}[rrdd]_2&&&\tquatretrois\ar@{-}[rd]^>(.4){1}&&
\ttroisun\tun\ar@{-}|(.50)\hole[llld]^>(.7){1}\ar@{-}[ld]^3\\
&&&\ttroisdeux\tun\ar@{-}[d]_3&&\tquatrequatre\ar@{-}[lld]^1&\\
&&&\tquatrecinq&&&}$$
The indices on the edges give the indices of the corresponding admissible transformation.
To obtain the Hasse graph of $\geq_I$, it is enough to delete the edges whose indices are not in $I$.

\subsection{Application to the product in the dual basis}

The coproduct of a forest can be expressed in terms of the order relations $\leq_I$:

\begin{prop}
Let $F \in \F-\{1\}$. Then:
$$\Delta(F)=F \otimes 1+1\otimes F+\sum_{\substack{F_1,F_2 \in \F-\{1\}\\F \leq_{\{weight(F_1)\}}F_1F_2}} F_1 \otimes F_2.$$
\end{prop}

{\bf Proof.} By induction on $n=weight(F)$. If $n=1$, then $F=\tun$ and the result is obvious.
Suppose the result true at ranks $\leq n-1$. Two cases are possible.
\begin{enumerate}
\item $F=F_1F_2$, and the induction hypothesis holds for $F_1$ and $F_2$. Then:
\begin{eqnarray*}
\tdelta(F)&=&(F_1 \otimes 1)\tdelta(F_2)+\tdelta(F_1)(1\otimes F_2)+F_1\otimes F_2\\
&=&\sum_{\substack{G_1,G_2 \in \F-\{1\}\\F_2 \leq_{\{weight(G_1)\}}G_1G_2}}F_1G_1 \otimes G_2
+\sum_{\substack{G_1,G_2 \in \F-\{1\}\\F_1 \leq_{\{weight(G_1)\}}G_1G_2}}G_1 \otimes F_2G_2+F_1 \otimes F_2\\
&=&\sum_{\substack{G_1,G_2 \in \F-\{1\}\\weight(G_1)>weight(F_1)\\F \leq_{\{weight(G_1)\}}G_1G_2}}G_1 \otimes G_2
+\sum_{\substack{G_1,G_2 \in \F-\{1\}\\weight(G_1)<weight(F_1)\\F \leq_{\{weight(G_1)\}}G_1G_2}}G_1 \otimes G_2\\
&&+\sum_{\substack{G_1,G_2 \in \F-\{1\}\\weight(G_1)=weight(F_1)\\F \leq_{\{weight(G_1)\}}G_1G_2}}G_1 \otimes G_2.
\end{eqnarray*}
So the result is true for $F$.
\item $F=B^+(G)$, and the induction hypothesis holds for $G$. Then:
\begin{eqnarray*}
\tdelta(F)&=&G\otimes \tun + (Id \otimes B^+)\circ \tdelta(G)\\
&=&G\otimes \tun +\sum_{\substack{F_1,F_2\in \F-\{1\}\\G \leq_{\{weight(F_1)\}} F_1F_2}} F_1 \otimes B^+(F_2)\\
&=&\sum_{\substack{F_1,F_2\in \F-\{1\}\\weight(F_2)=1\\ F \leq_{\{weight(F_1)\}} F_1F_2}} F_1 \otimes F_2
 +\sum_{\substack{F_1,F_2\in \F-\{1\}\\weight(F_2)>1\\ F \leq_{\{weight(F_1)\}} F_1F_2}} F_1 \otimes F_2.
\end{eqnarray*}
So the result is true for $F$. $\Box$ \\
\end{enumerate}

Dually:
\begin{cor}
\label{29}
Let $F_1,F_2 \in \F-\{1\}$. Then $\displaystyle f_{F_2}f_{F_1}=\sum_{\substack{G \in \F\\G \leq_{\{weight(F_1)\}} F_1F_2}} f_G$.
\end{cor}

{\bf Proof.} We put $\displaystyle f_{F_2}f_{F_1}=\sum_{G \in \F} a_G f_G$. Then:
\begin{eqnarray*}
a_G&=&\langle f_{F_2}f_{F_1},G \rangle\\
&=&\langle f_{F_1} \otimes f_{F_2},\Delta(G) \rangle\\
&=&\langle f_{F_1} \otimes f_{F_2},G \otimes 1+1\otimes G
+\sum_{\substack{G_1,G_2 \in \F-\{1\}\\G \leq_{\{weight(G_1)\}}G_1G_2}}G_1 \otimes G_2 \rangle\\
&=&\sum_{\substack{G_1,G_2 \in \F-\{1\}\\G \leq_{\{weight(G_1)\}}G_1G_2}}\delta_{F_1,G_1} \delta_{F_2,G_2}.
\end{eqnarray*}
So $a_G=1$ if $G \leq_{\{weight(F_1)\}} F_1F_2$, and $0$ if not. $\Box$\\

{\bf Examples.}
$$\left\{ \begin{array}{rcl}
f_{\tdeux} f_{\tun}&=&f_{\tun \tdeux}+f_{\ttroisun}+f_{\ttroisdeux},\\
f_{\tun}f_{\tdeux}&=&f_{\tdeux \tun}+f_{\ttroisdeux},\\
f_{\tdeux}f_{\tdeux}&=&f_{\tdeux\tdeux}+f_{\tquatredeux}+f_{\tquatrecinq}.
\end{array}\right.$$

The following result can be proved by induction on $n$:

\begin{cor}
Let $F_1,\ldots,F_n \in \F-\{1\}$. We define:
$$I=\{weight(F_1),weight(F_1F_2),\ldots,weight(F_1\ldots F_{n-1})\}.$$
Then:
$$f_{F_n} \ldots f_{F_1}=\sum_{\substack{G \in \F\\G \leq_I F_1\ldots F_n}}f_G.$$
\end{cor}

 For example, for $F_1=\ldots=F_n=\tun$, then $I=\{1,\ldots,n-1\}$, so $\geq_I=\geq$ on $\F(n)$.
By corollary \ref{23}, for all $F \in \F(n)$, $F \leq \tun^n$, so $\tun^n$ is the greatest element of $\F(n)$
(it is of course also possible to prove directly this result).

\subsection{Isomorphism of $(\F,\leq)$ with the Tamari poset}

Let $\T_b$ be the set of planar binary trees and, for all $n \in \mathbb{N}$, let $\T_b(n)$ be the set of planar binary trees with $n$ internal vertices
(and $n+1$ leaves). For example:
$$\T_b(0)=\{\bun\},\:\T_b(1)=\left\{\bdeux\right\},\:\T_b(2)=\left\{\btroisun,\btroisdeux\right\},\:
\T_b(3)=\left\{\bquatreun,\bquatredeux,\bquatretrois,\bquatrequatre,\bquatrecinq\right\}\ldots$$

The set $\T_b$  is given the Tamari order (see \cite{Aguiar,Ag2,Loday,Stanley2}).
Let us recall that it is the order generated by the local Tamari transformation:
$$\btroisun\longrightarrow \btroisdeux.$$

Let $t_1,t_2$ be two planar binary trees. We denote by $t_1\vee t_2$ the planar binary tree obtained by grafting $t_1$ and $t_2$ on a common root.
For example:
$$\bun \vee \btroisun=\bquatretrois,\: \bdeux \vee \bdeux=\bquatrecinq.$$
Note that every planar binary tree different from $\bun$ can be uniquely written $t^l \vee t^r$, where $t^l$ and $t^r$ are two planar binary trees. \\

We define a bijection $\eta:\T_b \longrightarrow \F$ by induction on the number of internal vertices by:
$$ \eta: \left\{\begin{array}{rcl}
\T_b& \longrightarrow &\F\\
\bun &\longrightarrow &1,\\
t_1 \vee t_2 & \longrightarrow &B^+(\eta(t_1)) \eta(t_2).
\end{array}\right. $$
It is not difficult to prove that $\eta$ is one-to-one, with inverse given by:
$$\eta^{-1}: \left\{\begin{array}{rcl}
\F& \longrightarrow & \T_b\\
1& \longrightarrow &\bun,\\
B^+(F_1)F_2&\longrightarrow& \eta^{-1}(F_1)\vee \eta^{-1}(F_2).
\end{array}\right. $$
Moreover, $\eta$ induces a bijection: $\eta:\T_b(n) \longrightarrow \F(n)$, for all $n \in \mathbb{N}$. \\

{\bf Examples.}
$$\begin{array}{rclc|crclc|crclc|crcl}
\bdeux& \stackrel{\eta}{\longrightarrow} &\tun&&&\btroisun&\stackrel{\eta}{\longrightarrow}&\tdeux&&&
\btroisdeux&\stackrel{\eta}{\longrightarrow}&\tun\tun&&&\bquatreun&\stackrel{\eta}{\longrightarrow}&\ttroisdeux\\
\bquatredeux&\stackrel{\eta}{\longrightarrow}&\ttroisun&&&\bquatretrois&\stackrel{\eta}{\longrightarrow}&\tun\tdeux&&&
\bquatrecinq&\stackrel{\eta}{\longrightarrow}&\tdeux\tun&&&\bquatrequatre&\stackrel{\eta}{\longrightarrow}&\tun\tun\tun
\end{array}$$

The aim of this subsection is to show the following theorem:
\begin{theo}
\label{31}
Let $t_1,t_2 \in \T_b$. The following equivalence holds:
$$(t_1 \leq t_2) \Longleftrightarrow (\eta(t_1)\leq \eta(t_2)).$$
In other terms, $\eta$ is a poset isomorphism.
\end{theo}

Let us define $\setminus:\T_b\times \T_b \longrightarrow \T_b$ in the following way (see \cite{Loday}):
$t_1 \setminus t_2$ is the grafting of  $t_2$ on the leave of $t_1$ at most on the right.
Note that $\setminus$ is an associative operation, and every $t \in \T_b$ can be uniquely decomposed in the form
$t=(t_1\vee \bun) \setminus \ldots \setminus (t_n \vee \bun)$, with $n \in \mathbb{N}$ and $t_1,\ldots,t_n \in \T_b$. \\

{\bf Examples.}
$$\bdeux \setminus \bdeux=\btroisdeux,\: \btroisun\setminus \bdeux =\bquatrecinq,\:\bdeux \setminus \btroisun= \bquatretrois.$$

\begin{lemma}
\label{32}
The Tamari order on $\T_b$ is the less fine partial order satisfying:
\begin{enumerate}
\item  for all $t_1,t_2,t_3 \in \T_b$, $(t_1 \leq t_2)\Longrightarrow \left\{\begin{array}{rcl}
t_1\vee \bun &\leq& t_2 \vee \bun,\\
t_1\setminus t_3 &\leq& t_2 \setminus t_3,\\
t_3 \setminus t_1 &\leq& t_3 \setminus t_2.
\end{array}\right. $
\item for all $t_1,t_2 \in \T_b$, $(t_1 \vee t_2)\vee \bun \leq t_1 \vee (t_2 \vee \bun)$.
\end{enumerate}
\end{lemma}

{\bf Proof.} Let $\sqsubseteq$ be the less fine partial order satisfying the assertions of lemma \ref{32}. 
Clearly, the Tamari order $\leq$ satisfies these conditions, so for all $t_1,t_2 \in \T_b$, $(t_1  \sqsubseteq t_2)\Longrightarrow (t_1 \leq t_2)$.

Let $t_1, t_2 \in \T_b$, such that $t_1 \leq t_2$. Let us show that $t_1  \sqsubseteq t_2$.
First, remark that $t_1$ and $t_2$ have the same number $n$ of internal vertices.
Let us proceed by induction on $n$. If $n=0$, then $t_1=t_2=\bun$ and the result is obvious.
If $n\geq 1$, we can suppose that $t_2$ is obtained from $t_1$ by a single Tamari operation $\btroisun \longrightarrow \btroisdeux$.
Let us put $t_1=t_1'\vee t_1''$ and $t_2=t_2'\vee t_2''$. Three cases are possible.
\begin{enumerate}
\item If the operation holds on a vertex of $t_1'$, then $t'_1 \leq t'_2$ and $t_1''=t_2''$. By induction hypothesis, $t_1'\leq t_2'$. 
So, $t_1'\vee \bun \sqsubseteq  t_2' \vee \bun$ and $t_1=(t_1'\vee \bun)\setminus t_1''  \sqsubseteq  (t_2' \vee \bun)\setminus t_1''=t_2$.
\item If the operation holds on a vertex of $t_1''$, then $t_1'=t_2'$ and $t_1'' \sqsubseteq t_2''$. 
By induction hypothesis, $t_1= (t_1' \vee \bun)\setminus t_1'' \sqsubseteq (t_2' \vee \bun)\setminus t_1''=t_2$.
\item If the operation holds on the root of $t_1$, let us put $t_1=(t_1'\vee t_1'')\vee t_1'''$. Then $t_2= t_1'\vee(t_1''\vee t_1''')$. 
Hence, $(t_1'\vee t_1'')\vee \bun \sqsubseteq   t_1'\vee(t_1''\vee \bun)$ and
$t_1=((t_1'\vee t_1'')\vee \bun)\setminus t_1''' \sqsubseteq  (t_1'\vee(t_1''\vee \bun))\setminus t_1'''=t_2$. $\Box$
\end{enumerate}

\begin{lemma}
\label{33}
The partial order $\leq$ on $\F$ is the less fine partial order satisfying:
\begin{enumerate}
\item  For all $F_1,F_2,F_3 \in \F$, $(F_1 \leq F_2)\Longrightarrow \left\{\begin{array}{rcl}
B^+(F_1) &\leq& B^+(F_2),\\
F_1F_3 &\leq& F_2F_3,\\
F_3F_1 &\leq& F_3 F_2.
\end{array}\right. $
\item For all $F_1,F_2 \in \F$, $B^+(B^+(F_1)F_2))\leq B^+(F_1)B^+(F_2)$.
\end{enumerate}
\end{lemma}

{\bf Proof.} Let $\sqsubseteq$ be the less fine partial order on $\F$ satisfying the assertions of lemma \ref{33}. 
As the order $\leq$ clearly satisfies these conditions, for all $F_1,F_2 \in \F$, $(F_1  \sqsubseteq F_2)\Longrightarrow (F_1 \leq F_2)$.

Let $F_1, F_2 \in \F$, such that $F_1 \leq F_2$. Let us show $F_1  \sqsubseteq F_2$. Necessarily, $F_1$ and $F_2$ have the same weight $n$.
Let us proceed by induction on $n$. If $n=0$, then $F_1=F_2=1$ and the result is obvious.
Suppose $n \geq 1$. We can suppose that $F_2$ is obtained from $F_1$ by a single elementary operation. Three cases are possible.
\begin{enumerate}
\item $F_1=B^+(G_1)$ and the operation holds on $G_1$. Then $F_2=B^+(G_2)$, with $G_1 \leq G_2$. By induction hypothesis, $G_1 \sqsubseteq  G_2$.
So, $F_1=B^+(G_1) \sqsubseteq  B^+(G_2)=F_2$.
\item $F_1=B^+(G_1)$ and the transformation is of second kind. We then put $G_1=t_1G_1'$ with $t_1 \in \T$. Then $F_2=t_1B^+(G_1')$.
By property  $2$ of $\sqsubseteq$, $F_1 \sqsubseteq F_2$.
\item $F_1=t_1\ldots t_k$, $k\geq 2$, and the transformation holds on $t_i$. 
We can then write $F_2=t_1\ldots t_{i-1}G_i t_{i+1}\ldots t_k$, with $t_i \leq G_i$. By induction hypothesis, $t_i \sqsubseteq G_i$. 
Hence, $t_1\ldots t_i \sqsubseteq t_1\ldots t_{i-1}G_i$ and $F_1=t_1\ldots t_k  \sqsubseteq  t_1\ldots t_{i-1}G_i t_{i+1}\ldots t_k=F_2$. $\Box$ 
\end{enumerate}

\begin{lemma}
\label{34}
$\eta$ satisfies the following assertions:
\begin{enumerate}
\item for all $t\in \T_b$, $\eta(t\vee \bun)=B^+(\eta(t))$.
\item for all $t_1,t_2\in \T_b$, $\eta(t_1\setminus t_2)=\eta(t_1)\eta(t_2)$.
\end{enumerate}
\end{lemma}

{\bf Proof.} \begin{enumerate}
\item Indeed, $\eta(t\vee \bun)=B^+(\eta(t)) \eta(\bun)=B^+(\eta(t))1=B^+(\eta(t))$.
\item Let us put $t_1= (s_1\vee \bun) \setminus \ldots \setminus (s_k\vee \bun)$. We proceed by induction on $k$.
If $k=0$, then $t_1=\bun$ and  $\eta(t_1)=1$, so $\eta(t_1\setminus t_2)=\eta(t_2)=\eta(t_1)\eta(t_2)$.
If $k=1$, then $t_1\setminus t_2=s_1 \vee t_2$ and $\eta(t_1)=B^+(\eta(s_1))$ by the first point. 
So, $\eta(t_1\setminus t_2)=B^+(\eta(s_1))\eta(t_2)=\eta(t_1)\eta(t_2)$. Let us suppose the result at rank $k-1$.
\begin{eqnarray*}
\eta(t_1\setminus t_2)&=& \eta((s_1\vee \bun) \setminus \ldots \setminus (s_k\vee \bun)\setminus t_2) \\
&=&\eta(s_1\vee\bun)  \eta((s_2\vee \bun) \setminus \ldots \setminus (s_k\vee \bun)\setminus t_2) \\
&=&\eta(s_1\vee\bun)  \eta(s_2\vee \bun) \ldots \eta(s_k\vee \bun)\eta(t_2) \\
&=&\eta((s_1\vee \bun) \setminus \ldots \setminus (s_k\vee \bun))\eta(t_2)\\
&=&\eta(t_1)\eta(t_2). 
\end{eqnarray*}
We used the result at rank $1$ for the second equality and the result at rank $k-1$ for the third and the fourth ones. $\Box$
\end{enumerate}

{\bf Proof of theorem \ref{31}.}
We define $\preceq$ on $\T_b$ by $(t_1\preceq t_2) \Longleftrightarrow (\eta(t_1)\leq \eta(t_2)).$ As $\eta$ is a bijection, this is a partial order on $\T_b$. 
By lemmas \ref{33} and \ref{34}, this is the less fine partial order satisfying the conditions of lemma \ref{32}. So it is the Tamari order. $\Box$

\subsection{A decreasing isomorphism of the poset $\F$}

\begin{prop}
We define an involution $m:\F\longrightarrow \F$ by induction on the weight in the following way:
$$ \left\{ \begin{array}{rcl}
m(1)&=&1,\\
m(B^+(F_1)F_2)&=&B^+(m(F_2))m(F_1)\mbox{ for any $F_1$, $F_2 \in \F $.}
\end{array}\right.$$
\end{prop}

{\bf Proof.} Clearly, this defines inductively $m(F)$ for all forest $F$ in a unique way.
Let us show that $m^2(F)=F$ for all forest $F \in \F$ by induction on the weight $n$ of $F$. If $n=0$, the result is obvious. 
Suppose the result true for all rank $<n$. Let $F \in \F(n)$. We put $F=B^+(F_1)F_2$. Then the induction hypothesis holds for $F_1$. So:
$$m\circ m(F)=m(B^+(m(F_2))m(F_1))=B^+(m\circ m(F_1))m\circ m(F_2)=B^+(F_1)F_2=F.$$
So $m$ is an involution. $\Box$\\

{\bf Remark.} For all forest $F \in \F$, $F$ and $m(F)$ have the same weight. So $m$ induces an involution 
$m:\F(n)\longrightarrow \F(n)$ for all $n \in \mathbb{N}$.\\

{\bf Examples.}
$$\begin{array}{rclc|crclc|crclc|crcl}
\tun&\stackrel{m}{\longleftrightarrow}&\tun&&&\tun\tun&\stackrel{m}{\longleftrightarrow}&\tdeux&&&
\tun\tun\tun&\stackrel{m}{\longleftrightarrow}&\ttroisdeux&&&\tun\tdeux&\stackrel{m}{\longleftrightarrow}&\ttroisun\\
\tdeux\tun&\stackrel{m}{\longrightarrow}&\tdeux\tun&&&\tun\tun\tun\tun&\stackrel{m}{\longleftrightarrow}&\tquatrecinq&&&
\tun\tun\tdeux&\stackrel{m}{\longleftrightarrow}&\tquatrequatre&&&\tun\tdeux\tun&\stackrel{m}{\longleftrightarrow}&\tquatredeux\\
\tun\ttroisun&\stackrel{m}{\longleftrightarrow}&\tquatretrois&&&\tun\ttroisdeux&\stackrel{m}{\longleftrightarrow}&\tquatreun&&&
\tdeux\tun\tun&\stackrel{m}{\longleftrightarrow}&\ttroisdeux\tun&&&\tdeux\tdeux&\stackrel{m}{\longleftrightarrow}&\ttroisun\tun
\end{array}$$

\begin{prop}
Let $F,G\in \F$. Then $F\leq G$ if, and only if, $m(G) \leq m(F)$.
\end{prop}

{\bf Proof.}
We consider the bijection $m':\T_b \longrightarrow \T_b$ defined by $m'=\eta^{-1} \circ m \circ \eta$. Then $m'(\bun)=\bun$, and for all $t_1,t_2 \in \T_b$:
\begin{eqnarray*}
m'(t_1 \vee t_2)&=&\eta^{-1} \circ m(B^+(\eta(t_1))\eta(t_2))\\
&=&\eta^{-1}(B^+(m\circ \eta(t_2))m\circ \eta(t_1))\\
&=&\eta^{-1}\circ m \circ \eta(t_2)\vee \eta^{-1} \circ m \circ \eta(t_1)\\
&=&m'(t_2) \vee m'(t_1).
\end{eqnarray*}
Hence, $m'$ is the vertical reflection:
$$\begin{array}{rclcrclcrclcrclcrcl}
\bdeux& \stackrel{m'}{\longleftrightarrow} &\bdeux,&\btroisun&\stackrel{m'}{\longleftrightarrow}&\btroisdeux,&
\bquatreun&\stackrel{m'}{\longleftrightarrow}&\bquatrequatre,&\bquatredeux&\stackrel{m'}{\longleftrightarrow}&\bquatretrois,&
\bquatrecinq&\stackrel{m'}{\longleftrightarrow}&\bquatrecinq
\end{array}$$
It is obviously a decreasing automorphism of $\T_b$. If $F,G \in \F$, let us put $t_1=\eta^{-1}(F)$ and $t_2=\eta^{-1}(G)$:
\begin{eqnarray*}
F \leq G & \Longleftrightarrow & \eta(t_1) \leq \eta(t_2) \\
& \Longleftrightarrow & t_1 \leq t_2 \\
& \Longleftrightarrow & m'(t_2) \leq m'(t_1) \\
& \Longleftrightarrow & \eta \circ m'(t_2) \leq \eta \circ m'(t_1) \\
& \Longleftrightarrow & \eta \circ m'\circ \eta^{-1}(G) \leq \eta \circ m'\circ \eta^{-1}(F)\\
& \Longleftrightarrow & m(G) \leq m(F).
\end{eqnarray*}
Hence, $m$ is decreasing. $\Box$\\

{\bf Remark.} So $m(\tun^n)=(B^+)^n(1)$ (ladder of weight $n$) is the smallest element of $\F(n)$.

\subsection{Link between the pairing and the partial order}

\begin{theo}
Let $F,G\in \F$. The following assertions are equivalent:
\begin{enumerate}
\item $\langle F,G\rangle \neq 0$.
\item $\langle F,G\rangle =1$.
\item $m(G) \leq F$.
\end{enumerate}
\end{theo}

{\bf Proof.}

$1 \Longrightarrow 3$. Let us suppose that $\langle F,G\rangle  \neq 0$ and let us show that $F \geq m(G)$.
We proceed by induction on $weight(F)=weight(G)=n$. If $n=0$, then $F=G=1$ and $F \geq m(G)$. We now suppose the result at all rank $<n$. 
We put $G=B^+(G_1)G_2$. Then the induction hypothesis holds for $G_1$ and $G_2$. Moreover, $m(G)=B^+(m(G_2))m(G_1)$.
Let $c$ be the unique left admissible cut of $F$ such that $P^c(F)$ and $G_2$ have the same weight. Then, by homogeneity of $\langle-,-\rangle$:
$$\langle F,G\rangle =\langle F,B(G_1)G_2\rangle =\langle \Delta(F),G_2\otimes B(G_1)\rangle=\langle P^c(F)\otimes R^c(F), G_2 \otimes B(G_1)\rangle=1.$$
As $\langle F,G\rangle \neq 0$, then $\langle P^c(F),G_2\rangle $  and $\langle R^c(F),B(G_1)\rangle $ are non zero.
By the induction hypothesis, $P^c(F) \geq m(G_2)$ and $R^c(F) \geq m(B(G_1))=\tun m(G_1)$.
As $c$ is a left admissible cut, we easily deduce that $F \geq B^+(m(G_2))m(G_1)=m(G)$.\\

$3 \Longrightarrow 2$. Let us suppose that $m(G) \leq F$ and let us show that $\langle F,G\rangle =1$.
We proceed by induction on the common weight $n$ of $G$ and $F$. If $n=0$, then $F=G=1$ and $\langle F,G\rangle =1$. 
We now suppose the result at all rank $<n$. We put $G=B^+(G_1)G_2$, with $G_1,G_2 \in \F$.
Then $m(G)=B^+(m(G_2))m(G_1)$ and $F\geq B^+(m(G_2)) m(G_1)$.
Let $c$ be the unique left admissible cut of $F$ such that $P^c(F)$ and $G_2$ have the same weight.
By definition of $\geq$, $P^c(F) \geq m(G_2)$ and $R^c(F)\geq \tun m(G_1)=m(B^+(G_1))$.
Hence, $\langle P^c(F),G_2\rangle =\langle R^c(F),B^+(G_1)\rangle =1$. As $\langle-,-\rangle $ is homogeneous:
$$\langle F,G\rangle =\langle F,B^+(G_1)G_2\rangle =\langle \Delta_0(F),G_2\otimes B^+(G_1)\rangle
=\langle P^c(F)\otimes R^c(F), G_2 \otimes B^+(G_1)\rangle=1.$$

$2\Longrightarrow 1$. Obvious. $\Box$\\

{\bf Examples.} Values of the pairing $\langle-,-\rangle $ for forests of weight $\leq 4$:
$$\begin{array}{c|c}
&\tun\\
\hline
\tun&1
\end{array}\hspace{1cm}
\begin{array}{c|cc}
&\tun\tun&\tdeux\\
\hline
\tun\tun&1&1\\
\tdeux&1&0
\end{array}\hspace{1cm}
\begin{array}{c|ccccc}
&\tun\tun\tun&\tun\tdeux&\tdeux\tun&\ttroisun&\ttroisdeux\\
\hline
\tun\tun\tun&1&1&1&1&1\\
\tun\tdeux&1&1&0&1&0\\
\tdeux\tun&1&0&1&0&0\\
\ttroisun&1&1&0&0&0\\
\ttroisdeux&1&0&0&0&0
\end{array}$$
$$\begin{array}{c|cccccccccccccc}
&\tun\tun\tun\tun&\tun\tun\tdeux&\tun\tdeux\tun&\tun\ttroisun&
\tun\ttroisdeux&\tdeux\tun\tun&\tdeux\tdeux&\ttroisun\tun&\ttroisdeux\tun&
\tquatreun&\tquatretrois&\tquatredeux&\tquatrequatre&\tquatrecinq\\
\hline
\tun\tun\tun\tun&1&1&1&1&1&1&1&1&1&1&1&1&1&1\\
\tun\tun\tdeux&1&1&1&1&1&0&0&1&0&1&1&0&1&0\\
\tun\tdeux\tun&1&1&0&1&0&1&1&0&0&1&0&1&0&0\\
\tun\ttroisun&1&1&1&1&1&0&0&0&0&1&1&0&0&0\\
\tun\ttroisdeux&1&1&0&1&0&0&0&0&0&1&0&0&0&0\\
\tdeux\tun\tun&1&0&1&0&0&1&0&1&1&0&0&0&0&0\\
\tdeux\tdeux&1&0&1&0&0&0&0&1&0&0&0&0&0&0\\
\ttroisun\tun&1&1&0&0&0&1&1&0&0&0&0&0&0&0\\
\ttroisdeux\tun&1&0&0&0&0&1&0&0&0&0&0&0&0&0\\
\tquatreun&1&1&1&1&1&0&0&0&0&0&0&0&0&0\\
\tquatredeux&1&1&0&1&0&0&0&0&0&0&0&0&0&0\\
\tquatretrois&1&0&1&0&0&0&0&0&0&0&0&0&0&0\\
\tquatrequatre&1&1&0&0&0&0&0&0&0&0&0&0&0&0\\
\tquatrecinq&1&0&0&0&0&0&0&0&0&0&0&0&0&0
\end{array}$$\\

As a corollary, we can give another (shorter) proof of the symmetry of $\langle-,-\rangle $: if $F,G \in \F$,
$$(\langle F,G\rangle =1)\Longleftrightarrow (F\geq m(G)) \Longleftrightarrow (m^2(F)\geq m(G)) \Longleftrightarrow 
(G \geq m(F)) \Longleftrightarrow (\langle G,F\rangle =1).$$

\begin{cor}
Let $F\in \F(n)$. Then $\displaystyle F= \sum_{\substack{G\in \F(n)\\ F\geq m(G)}} f_G= \sum_{\substack{G\in \F(n)\\ F\geq G}} f_{m(G)}$.
\end{cor}

{\bf Proof.} We put $\displaystyle F=\sum_{G\in \F(n)} a_{G,F} f_G$. Then $a_{G,F}=\langle F,G\rangle $, which implies this corollary. $\Box$ \\

Let $\mu:\F(n)^2\longrightarrow K$ be the M\"obius function of the poset $\F(n)$, that is to say (see \cite{Stanley}):
\begin{enumerate}
\item $\mu(F,G)=0$ if $F\nleq G$.
\item $\displaystyle \sum_{\substack{G\in \F(n)\\ F\leq G\leq H}} \mu(F,G)=\delta_{F,H}$ if $F\leq H$.
\end{enumerate}
Immediately, by \cite{Stanley}:

\begin{cor}
\label{39}
Let $F\in \F(n)$. Then $\displaystyle f_F= \sum_{\substack{G\in \F(n)\\ G\leq m(F)}} \mu(G,m(F)) G$.
\end{cor}

{\bf Examples.}
$$\begin{array}{|c|}
\hline
\begin{array}{rclc|rcl}
f_1&=&1&&
f_{\tun}&=&\tun\\
f_{\tun\tun}&=&\tdeux&&
f_{\tdeux}&=&-\tdeux+\tun\tun\\
f_{\tun\tun\tun}&=&\ttroisdeux&&
f_{\tun\tdeux}&=&-\ttroisdeux+\ttroisun\\
f_{\tdeux\tun}&=&-\ttroisdeux+\tdeux\tun&&
f_{\ttroisun}&=&-\ttroisun+\tun\tdeux\\
f_{\ttroisdeux}&=&\ttroisdeux-\tdeux\tun-\tun\tdeux+\tun\tun\tun&&
f_{\tun\tun\tun\tun}&=&\tquatrecinq\\
f_{\tun\tun\tdeux}&=&-\tquatrecinq+\tquatrequatre&&
f_{\tun\tdeux\tun}&=&-\tquatrecinq+\tquatredeux\\
f_{\tun\ttroisun}&=&-\tquatrequatre+\tquatretrois&&
f_{\tun\ttroisdeux}&=&\tquatrecinq-\tquatredeux-\tquatretrois+\tquatreun\\
f_{\tdeux\tun\tun}&=&-\tquatrecinq+\ttroisdeux\tun&&
f_{\tdeux\tdeux}&=&\tquatrecinq-\tquatrequatre-\ttroisdeux\tun+\ttroisun\tun\\
f_{\ttroisun\tun}&=&-\tquatredeux+\tdeux\tdeux&&
f_{\ttroisdeux\tun}&=&\tquatrecinq-\ttroisdeux\tun-\tdeux\tdeux+\tdeux\tun\tun\\
f_{\tquatreun}&=&-\tquatretrois+\tun\ttroisdeux&&
f_{\tquatretrois}&=&\tquatretrois-\tquatreun-\tun\ttroisdeux+\tun\ttroisun\\
f_{\tquatredeux}&=&\tquatrequatre-\ttroisun\tun-\tun\ttroisdeux+\tun\tdeux\tun&&
f_{\tquatrequatre}&=&\tquatredeux-\tdeux\tdeux-\tun\ttroisun+\tun\tun\tdeux\\
\end{array}\\ \hline
f_{\tquatrecinq}=-\tquatrecinq+\ttroisdeux\tun+\tdeux\tdeux-\tdeux\tun\tun+
\tun\ttroisdeux-\tun\tdeux\tun-\tun\tun\tdeux+\tun\tun\tun\tun.\\
\hline \end{array}$$

\bibliographystyle{amsplain}
\bibliography{biblio}

\providecommand{\bysame}{\leavevmode\hbox to3em{\hrulefill}\thinspace}
\providecommand{\MR}{\relax\ifhmode\unskip\space\fi MR }
\providecommand{\MRhref}[2]{%
  \href{http://www.ams.org/mathscinet-getitem?mr=#1}{#2}
}
\providecommand{\href}[2]{#2}
\begin{thebibliography}{10}

\bibitem{Abe}
Eiichi Abe, \emph{Hopf algebras}, Cambridge Tracts in Mathematics, no.~74,
  Cambridge University Press, Cambridge-New York, 1980.

\bibitem{Aguiar3}
Marcelo Aguiar, \emph{Infinitesimal {H}opf algebras}, Contemp. Math.
  \textbf{267} (2000), 1--29.

\bibitem{Aguiar}
Marcelo Aguiar and Frank Sottile, \emph{Structure of the
  {M}alvenuto-{R}eutenauer {H}opf algebra of permutations}, Adv. Math.
  \textbf{191} (2005), no.~2, 225--275, math.CO/02 03282.

\bibitem{Ag2}
\bysame, \emph{Structure of the {L}oday-{R}onco {H}opf algebra of trees}, J.
  Algebra \textbf{295} (2006), no.~2, 473--511, math.CO/04 09022.

\bibitem{Connes}
Alain Connes and Dirk Kreimer, \emph{Hopf algebras, {R}enormalization and
  {N}oncommutative geometry}, Comm. Math. Phys \textbf{199} (1998), no.~1,
  203--242, hep-th/98 08042.

\bibitem{Connes2}
\bysame, \emph{Renormalization in quantum field theory and the
  {R}iemann-{H}ilbert problem {I}. {T}he {H}opf algebra of graphs and the main
  theorem}, Comm. Math. Phys. \textbf{210} (2000), no.~1, 249--273, hep-th/99
  12092.

\bibitem{Connes3}
\bysame, \emph{Renormalization in quantum field theory and the
  {R}iemann-{H}ilbert problem. {II}. {T}he $\beta$-function, diffeomorphisms
  and the renormalization group}, Comm. Math. Phys. \textbf{216} (2001), no.~1,
  215--241, hep-th/00 03188.

\bibitem{Foissy}
Lo{\"\i}c Foissy, \emph{Finite-dimensional comodules over the {H}opf algebra of
  rooted trees}, J. Algebra \textbf{255} (2002), no.~1, 85--120, math.QA/01
  05210.

\bibitem{Foissy2}
\bysame, \emph{Les alg\`ebres de {H}opf des arbres enracin\'es, {I}}, Bull.
  Sci. Math. \textbf{126} (2002), 193--239.

\bibitem{Holtkamp}
Ralf Holtkamp, \emph{Comparison of {H}opf {A}lgebras on {T}rees}, Arch. Math.
  (Basel) \textbf{80} (2003), no.~4, 368--383.

\bibitem{Kreimer1}
Dirk Kreimer, \emph{On the {H}opf algebra structure of pertubative quantum
  field theories}, Adv. Theor. Math. Phys. \textbf{2} (1998), no.~2, 303--334,
  q-alg/97 07029.

\bibitem{Kreimer2}
\bysame, \emph{On {O}verlapping {D}ivergences}, Comm. Math. Phys. \textbf{204}
  (1999), no.~3, 669--689, hep-th/98 10022.

\bibitem{Kreimer3}
\bysame, \emph{Combinatorics of (pertubative) {Q}uantum {F}ield {T}heory},
  Phys. Rep. \textbf{4--6} (2002), 387--424, hep-th/00 10059.

\bibitem{Loday}
Jean-Louis Loday and Maria~O. Ronco, \emph{Hopf algebra of the planar binary
  trees}, Adv. Math. \textbf{139} (1998), no.~2, 293--309.

\bibitem{Loday2}
\bysame, \emph{On the structure of cofree hopf algebras}, J. Reine Angew. Math.
  \textbf{592} (2006), 123--155.

\bibitem{Moerdijk}
Ieke Moerdijk, \emph{On the {C}onnes-{K}reimer construction of {H}opf
  algebras}, Contemp. Math. \textbf{271} (2001), 311--321, math-ph/99 07010.

\bibitem{Stanley}
Richard~P. Stanley, \emph{Enumerative combinatorics. {V}ol. 1.}, Cambridge
  Studies in Advanced Mathematics, no.~49, Cambridge University Press,
  Cambridge, 1997.

\bibitem{Stanley2}
\bysame, \emph{Enumerative combinatorics. {V}ol. 2.}, Cambridge Studies in
  Advanced Mathematics, no.~62, Cambridge University Press, Cambridge, 1999.

\bibitem{Sweedler}
Moss~E. Sweedler, \emph{Hopf algebras}, Mathematics Lecture Note Series, W. A.
  Benjamin, Inc., New York, 1969.

\end{thebibliography}

\end{document}